\documentclass[xypic]{amsart} 
\usepackage[latin1]{inputenc}
\usepackage[english]{babel}
\usepackage{amsfonts,amssymb,amsmath,latexsym,stmaryrd,upref}
\usepackage[active]{srcltx}
\usepackage[margin=2.5cm]{geometry}
\usepackage{graphicx}

\input xypic
\xyoption{all}
\CompileMatrices

\usepackage{color}
\newcommand{\txtr}{\textcolor{red}}

\newcommand{\txtb}{\textcolor{blue}}

\newtheorem{thm}{\txtb{Theorem}}[section] 
\newtheorem{lemma}[thm]{\txtb{Lemma}}
\newtheorem{cor}[thm]{\txtb{Corollary}}
\newtheorem{prop}[thm]{\txtb{Proposition}}

\theoremstyle{definition}
\newtheorem{defn}[thm]{\txtr{Definition}}
\newtheorem{exmp}[thm]{\txtb{Example}}
\newtheorem{exc}[thm]{\txtb{Exercise}}
\newtheorem{rmk}[thm]{\txtb{Remark}}

\newcommand{\func}[3]{\mbox{$#1\colon #2\rightarrow #3$}}

\newcommand{\set}{\setcounter{equation}{\value{thm}}}
\newcommand{\add}{\addtocounter{thm}{1}}

\newcommand{\n}{neighborhood}
\newcommand{\co}{continuous}

\renewcommand{\int}{\operatorname{int}}

\newcommand{\m}{morphism}

\newcommand{\Hom}{\operatorname{Hom}}

\newcommand{\U}{\mathrm{U}}
\newcommand{\SU}{\mathrm{SU}}
\newcommand{\SO}{\mathrm{SO}}
\newcommand{\PSO}{\mathrm{PSO}}
\newcommand{\Or}{\mathrm{O}}

\newcommand{\Spin}{\mathrm{Spin}}
\newcommand{\Symp}{\mathrm{Sp}}

\newcommand{\gen}[1]{\left\langle {#1} \right\rangle}

\newcommand{\order}[1]{\lvert #1 \rvert}
\newcommand{\diag}{\operatorname{diag}}
\newcommand{\lcm}{\operatorname{lcm}}

\newcommand{\C}{{\mathbf C}}
\newcommand{\Ha}{{\mathbf H}}
\newcommand{\Z}{{\mathbf Z}}

\newcommand{\R}{{\mathbf R}}

\newcommand{\GY}{G \backslash Y}

\numberwithin{equation}{section} 

\usepackage{floatflt}
\usepackage{color}

\usepackage[bookmarks=true,bookmarksopen=false]{hyperref}
\hypersetup{
   pdftitle = {},
   pdfauthor = {Jesper Michael Møller},
   pdfpagemode = {UseOutlines},
   pdfstartview = {FitH},
   pdfborder = {0 0 0},
   backref = {true},
   colorlinks = {true},
   urlcolor = {blue},
   citecolor = {blue},
   linkcolor = {blue},
   pdftoolbar = {true},
}

\usepackage{graphicx}

\begin{document}
\title{The fundamental group and covering spaces}

\author[J.M. Møller]{Jesper M.\ Møller}
\address{Matematisk Institut\\
  Universitetsparken 5\\
  DK--2100 København} 
\email{moller@math.ku.dk} 
\urladdr{\url{http://www.math.ku.dk/~moller}} 
\date{\today}
\begin{abstract}
  These notes, from a first course in algebraic topology, introduce
  the fundamental group and the fundamental groupoid of a topological
  space and use them to classify covering spaces.
\end{abstract}
\thanks{I would like to thank Morten Poulsen for supplying the
  graphics and Yaokun Wu for spotting some errors in an earlier
  version.  The author is partially supported by the DNRF through the
  Centre for Symmetry and Deformation in Copenhagen}
\maketitle
\tableofcontents

\section{Homotopy theory of paths and loops}
\label{sec:pathhtpy}

\begin{defn}
A {\em path\/} in a topological space $X$ from $x_0\in X$ to $x_1 \in
X$ is a map \func{u}{I}{X} of the unit interval into $X$ with
$u(0)=x_0$ and $u(1)=x_1$.

Two paths, $u_0$ and $u_1$, from $x_0$ to $x_1$ are {\em path
  homotopic}, and we write simply $u_0 \simeq u_1$, if $u_0 \simeq
u_1\ \mathrm{rel\ } \partial I$, ie if $u_0$ and $u_1$ are are homotopic
relative to the end-points $\{0,1\}$ of the unit interval $I$.  
\end{defn}
 
\begin{itemize}
\item The {\em constant path\/} at $x_0$ is the
path $x_0(s)=x_0$ for all $s \in I$
\item The {\em inverse path\/} to $u$
is the path from $x_1$ to $x_0$ given by $\overline{u}(s)=u(1-s)$ 
\end{itemize}
If $v$ is a path from $v(0)=u(1)$ then the {\em product path\/}
path $u \cdot v$ given by
\begin{equation*}
  (u \cdot v)(s)= 
  \begin{cases}
    u(2s) & 0 \leq s \leq \frac{1}{2}\\
    v(2s-1) &  \frac{1}{2} \leq s \leq 1\\
  \end{cases}
\end{equation*}
where we first run along $u$ with double speed and then along $v$ with
double speed is a path from $u(0)$ to $v(1)$.

In greater detail, $u_0 \simeq u_1$ if there exists a homotopy \func
h{I \times I}X such that $h(s,0)=u_0(s)$, $h(s,1)=u_1(s)$ and
$h(0,t)=x_0$, $h(1,t)=x_1$ for all $s,t \in I$.  All paths in a
homotopy class have the same start point and the same end point.
\begin{figure}[t]
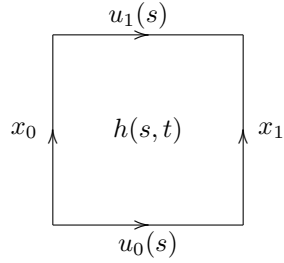

  \centering
\begin{equation*}
  \xy 0;/r.15pc/: 
  (-20,-20)*{}="A"; (20,-20)*{}="B"; (-20,20)*{}="C"; (20,20)*{}="D";
   "A"; "B" **\dir{-} 
   \POS?(.5)*{}*\dir{>}+(0,-4)*{u_0(s)};
   "C"; "D" **\dir{-}?(.5)*\dir{>}+(-2,4)*{u_1(s)}; 
   "A"; "C" **\dir{-}?(.5)*\dir{>}+(-6,0)*{x_0};
   "B"; "D" **\dir{-}?(.5)*\dir{>}+(6,0)*{x_1};
   (0,0)*{h(s,t)};
   \endxy
\end{equation*}  
  \caption{A path homototopy between two paths}
  \label{fig:pathhtpy}
\end{figure}
Note the following rules for products of paths
\begin{itemize}
\item $x_0 \cdot u \simeq u \simeq u \cdot x_1$
\item $u \cdot \overline{u} \simeq x_0$, $\overline{u} \cdot u \simeq x_1$
\item $(u \cdot v) \cdot w \simeq u \cdot (v \cdot w)$
\item $u_0 \simeq u_1,\ v_0 \simeq v_1 \Longrightarrow u_0
\cdot v_0 \simeq u_1\cdot v_1$
\end{itemize}
These drawings are meant to suggest proofs for the first three statements:
\begin{equation*}
  \xy 0;/r.15pc/: 
  (-20,-20)*{}="A"; (20,-20)*{}="B"; (-20,20)*{}="C"; (20,20)*{}="D";
   "A"; "B" **\dir{-} 
   \POS?(.5)*{}="-50" 
   \POS?(.25)*{}*\dir{>}+(0,-4)*{u(2s)}
   \POS?(.75)*{}*\dir{>}+(0,-4)*{x_1};
   "C"; "D" **\dir{-}?(.5)*\dir{>}+(-2,4)*{u(s)}; 
   "A"; "C" **\dir{-}?(.5)*\dir{>}+(-6,0)*{x_0};
   "B"; "D" **\dir{-}?(.5)*\dir{>}+(6,0)*{x_1};
   "-50"; "D" **\dir{.}; 
   \endxy \qquad
   \xy 0;/r.15pc/: 
   (-20,-20)*{}="A"; (20,-20)*{}="B"; (-20,20)*{}="C"; (20,20)*{}="D";
   "A"; "B" **\dir{-} 
   \POS?(.5)*{}="-50" 
   \POS?(.25)*{}*\dir{>}+(0,-4)*{u(2s)}
   \POS?(.75)*{}*\dir{>}+(0,-4)*{u(1-2s)};
   "C"; "D" **\dir{-}?(.5)*\dir{>}+(-2,4)*{x_0}; 
   "A"; "C" **\dir{-}?(.5)*\dir{>}+(-6,0)*{x_0};
   "B"; "D" **\dir{-}?(.5)*\dir{>}+(6,0)*{x_0};
   "-50"; "C" **\dir{.};
   "-50"; "D" **\dir{.};
   \endxy \qquad
   \xy 0;/r.15pc/: 
   (-20,-20)*{}="A"; (20,-20)*{}="B"; (-20,20)*{}="C"; (20,20)*{}="D";
   "A"; "B" **\dir{-} 
   ?(.25)+(-4,-4)*{u}?(.50)+(-4,-4)*{v}?(.75)+(-2,-4)*{w}
   \POS?(.25)*{}="-25" 
   \POS?(.50)*{}="-50";
   "C"; "D" **\dir{-}
   ?(.25)+(-2,4)*{u}?(.50)+(4,4)*{v}?(.75)+(4,4)*{w}
   \POS?(.50)*{}="+50" 
   \POS?(.75)*{}="+75";
   "A"; "C" **\dir{-}?(.5)*\dir{>}+(-6,0)*{x_0};
   "B"; "D" **\dir{-}?(.5)*\dir{>}+(6,0)*{x_3};
   "-25"; "+50" **\dir{.}?(.5)+(-4,0)*{x_1};
   "-50"; "+75" **\dir{.}?(.5)+(4,0)*{x_2};
   \endxy 
\end{equation*}

In the first case, we first run along $u$ with double speed and then
stand still at the end point $x_1$ for half the time. The homotopy consists
in slowing down on the path $u$ and spending less time just standing
still at $x_1$.

In the second case, we first run along $u$ all the way to $x_1$ with
double speed and then back again along $u$ also with double speed. The
homotopy consists in running out along $u$, standing still for an
increasing length of time (namely for $\frac{1}{2}(1-t) \leq s \leq
\frac{1}{2}(1+t)$), and running back along $u$.

In the third case we first run along $u$ with speed $4$,
then along $v$ with speed $4$, then along $w$ with speed $2$, and
we must show that can deformed into the case where we run along $u$
with speed $2$, along $v$ with speed $4$, along $w$ with speed
$4$. This can be achieved by slowing down on $u$, keeping the same
speed along $v$, and speeding up on $w$.

The fourth of the above rules is proved by this picture,
\begin{equation*}
   \xy 0;/r.15pc/: 
  (-20,-20)*{}="A"; (20,-20)*{}="B"; (-20,20)*{}="C"; (20,20)*{}="D";
   (60,-20)*{}="E"; (60,20)*{}="F";
   "A"; "E" **\dir{-}
   \POS?(.25)*\dir{>}+(0,-4)*{u_0(2s)}
   \POS?(.75)*\dir{>}+(0,-4)*{v_0(2s-1)};
   "C"; "F" **\dir{-}
   \POS?(.25)*\dir{>}+(0,4)*{u_1(2s)}
   \POS?(.75)*\dir{>}+(0,4)*{v_1(2s-1)};
   "A"; "C" **\dir{-}
   \POS?(.5)*\dir{>}+(-4,0)*{x_0};  
   "B"; "D" **\dir{-}
   \POS?(.5)*\dir{>}+(-4,0)*{x_1}; 
   "E"; "F" **\dir{-}
   \POS?(.5)*\dir{>}+(-4,0)*{x_2}; 
   \endxy
\end{equation*}

We write $\pi(X)(x_0,x_1)$ for the set of al homotopy classes of paths
in $X$ from $x_0$ to $x_1$ and we write $[u] \in \pi(X)(x_0,x_1)$ for
the homotopy class containing the path $u$.  The fourth rule implies
that the product operation on paths induces a product operation on
homotopy classes of paths
\set\begin{equation}\label{eq:defpathprod}
  \pi(X)(x_0,x_1) \times \pi(X)(x_1,x_2) \xrightarrow{\cdot} 
  \pi(X)(x_0,x_2) \colon ([u],[v]) \to [u] \cdot [v]= [u \cdot v]
\end{equation}\add
and the other three rules translate to similar rules
\begin{itemize}
\item $[x_0] \cdot [u] = [u] = [u] \cdot [x_1]$ (neutral elements)
\item $[u] \cdot [\overline{u}] = [x_0]$, $[\overline{u}] \cdot [u] =
  [x_1]$ (inverse elements)
\item $([u] \cdot [v]) \cdot [w] = [u] \cdot ([v] \cdot [w])$
  (associativity) 
\end{itemize}
for this product operation.

We next look at the functorial properties of this construction. 
Suppose that \func{f}{X}{Y} is a map of spaces. If $u$ is a path in
$X$ from $x_0$ to $x_1$, then the image $fu$ is a path in $Y$ from
$f(x_0)$ to $f(x_1)$. Since homotopic paths have homotopic images
there is an induced map
\begin{equation*}
  \pi(X)(x_0,x_1) \xrightarrow{\pi(f)} \pi(Y)(f(x_0),f(x_1)) \colon
  [u] \to [fu]
\end{equation*}
on the set of homotopy classes of paths. Observe that this map
\begin{itemize}
\item  $\pi(f)$  does not change if we change $f$ by a homotopy
relative to $\{x_0,x_1\}$,
\item  $\pi(f)$ respects the
path product operation in the sense that $\pi(f)([u] \cdot
[v])=\pi(f)([u]) \cdot \pi(f)([v])$ when $u(1)=v(0)$,
\item $\pi(\mathrm{id}_X)=\mathrm{id}_{\pi(X)(x_1,x_2)}$, $\pi(g \circ
  f) = \pi(g) \circ \pi(f)$ for maps \func gYZ
\end{itemize}

We now summarize our findings.

\begin{prop}
  For any space $X$, $\pi(X)$ is a groupoid, and for any map \func fXY
  between spaces the induced map \func{\pi(f)}{\pi(X)}{\pi(Y)} is a
  groupoid homo\m . In fact, $\pi$ is a functor from the category of
  topological spaces to the category of groupoids.
\end{prop}

\begin{defn} 
  $\pi(X)$ is called the fundamental groupid of $X$.  The {\em
    fundamental group\/} based at $x_0 \in X$ is the group
  $\pi_1(X,x_0)=\pi(X)(x_0,x_0)$ of homotopy classes of loops in $X$
  based at $x_0$.
\end{defn}

The path product \eqref{eq:defpathprod} specializes to a product
operation
\begin{equation*}
  \pi_1(X,x_0) \times \pi_1(X,x_0) \to \pi_1(X,x_0)
\end{equation*}
and to
transitive free group actions
\set\begin{equation}\label{eq:pi1actions}
  \pi_1(X,x_0) \times \pi(X)(x_0,x_1) \to  \pi(X)(x_0,x_1) \leftarrow
   \pi(X)(x_0,x_1) \times \pi_1(X,x_1) 
\end{equation}\add
so that $\pi_1(X,x_0)$ is indeed a group and $\pi(X)(x_0,x_1)$ is an
affine group from the left and from the right.

For fundamental groups, in particular, any based map
\func{f}{(X,x_0)}{(Y,y_0)} induces a group homo\m\
\func{\pi(f)=f_*}{\pi_1(X,x_0)}{\pi_1(Y,y_0)}, given by
$\pi_1(f)=f_*([u])=[fu]$, that only depends on the based homotopy
class of the based map $f$.

\begin{prop}
  The fundamental group is a functor
  \func{\pi_1}{\mathbf{hoTop}_*}{\mathbf{Grp}} from the homotopy
  category of based topological spaces into the category of groups.
\end{prop}

This means that $\pi_1(\mathrm{id}_{(X,x_0)}) =
\mathrm{id}_{\pi_1(X,x_0)}$ and $\pi_1(g \circ f) = \pi_1(g) \circ
\pi_1(f)$.  It follows immediately that if \func fXY is a homotopy
equivalence of {\em based\/} spaces then the induced map \func
{\pi_1(f)=f_*}{\pi_1(X,x_0)}{\pi_1(Y,y_0)} is an iso\m\ of groups.
(See Section~\ref{sec:cat} for more information about categories and
functors.)

\begin{cor}\label{cor:piretract}
  Let $X$ be a space, $A$ a subspace, and
  \func{i_*}{\pi_1(A,a_0)}{\pi_1(X,a_0)} the group homo\m\ induced by
  the inclusion map \func iAX.
  \begin{enumerate}
  \item If $A$ is a retract of $X$ then $i_*$ has a left inverse (so
    it is a mono\m ).
  \item If $A$ is a deformation retract of $X$ then $i_*$ has an
    inverse (so it is an iso\m ).
  \end{enumerate}
\end{cor}
\begin{proof}
  (1) Let \func rXA be a map such that $ri=1_A$. Then $r_*i_*$ is the
  identity iso\m\ of $\pi_1(A,a_0)$.

  \noindent (2) Let \func rXA be a map such that $ri=1_A$ and $ir
  \simeq 1_X \mathrm{rel\ }A$. Then
  $r_*i_*$ is the identity iso\m\ of $\pi_1(A,a_0)$  and $i_*r_*$ is
  the identity iso\m\ of $\pi_1(X,a_0)$ so $i_*$ and $r_*$ are each
  others' inverses.
\end{proof}

\begin{cor}
  Let $X$ and $Y$ be spaces. There is an iso\m\
  \begin{equation*}
    \func {(p_X)_* \times (p_Y)_*}{\pi_1(X \times Y,x_0 \times y_0)}
    {\pi_1(X,x_0) \times \pi_1(Y,y_0)}
  \end{equation*}
  induced by the projection maps \func {p_X}{X \times Y}X and \func
  {p_Y}{X \times Y}Y.  
\end{cor}
\begin{proof}
  The loops in $X \times Y$ have the form $u \times v$ where $u$ and
  $v$ are loops in $X$ and $Y$, respectively
  (\href{http://www.math.ku.dk/~moller/e03/3gt/notes/gtnotes.pdf}{General
    Topology, 2.63}). 
  The above homo\m\ has the form $[u \times v] \to
  [u] \times [v]$. The inverse homo\m\ is $[u] \times [v] \to [u
  \times v]$. Note that this is well-defined.
\end{proof}

We can now compute our first fundamental group.

\begin{exmp}\label{exmp:pi1Rn}
  $\pi_1(\R^n,0)$ is the trivial group with just one element because
  $\R^n$ contains the subspace $\{0\}$ consisting of one point as a
  deformation retract. Any space that deformation retract onto on of
  its points has trivial fundamental group. Is it true that any
  contractible space has trivial fundamental group?
\end{exmp}

Our tools to compute $\pi_1$ in more interesting cases are covering
space theory and van Kampen's theorem.

\setcounter{subsection}{\value{thm}}
\subsection{Change of base point and unbased homotopies}
\label{sec:unbased}\addtocounter{thm}{1}

What happens if we change the base point? In case, the new base point
lies in another path-component of $X$, there is no relation at all
between the fundamental groups. But if the two base points lie in the
same path-component, the fundamental groups are isomorphic.

\begin{lemma}
  If  $u$ is a path from  $x_0$ to $x_1$ then conjugation with $[u]$ 
  \begin{equation*}
    \pi_1(X,x_1) \to \pi_1(X,x_0) \colon [v] \to [u] \cdot [v] \cdot
    [\overline{u}] 
  \end{equation*}
  is a group iso\m .
\end{lemma}
\begin{proof}
  This is immediate from the rules for products of paths and a special
  case of \eqref{eq:pi1actions}.
\end{proof}

We already noted that if two maps are homotopic \emph{relative to the
  base point\/} then they induce the same group homo\m\ between the
fundamental groups.  We shall now investigate how the fundamental
group behaves with respect to free maps and free homotopies, ie maps
and homotopies that do not preserve the base point.

\begin{lemma}\label{lemma:unbasedmaps}
   Suppose
  that
  \func{f_0 \simeq f_1}XY are homotopic maps and \func h{X \times I}Y
  a homotopy. For any point $x \in X$, let $h(x) \in
  \pi(Y)(f_0(x),f_1(x))$ be the path homotopy class of $t \to
  h(x,t)$. For any $u \in \pi(X)(x_0,x_1)$
  there is a commutative diagram 
  \begin{equation*}
    \xymatrix{f_0(x_0) \ar[d]_{f_0(u)} \ar[r]^{h(x_0)} &
      f_1(x_0) \ar[d]^{f_1(u)} \\
     f_0(x_1) \ar[r]_{h(x_1)} & f_1(x_1) }
  \end{equation*}
  in $\pi(Y)$.
\end{lemma}
\begin{proof}
  Let $u$ be any path from $x_0$ to $x_1$ in $X$.  If we push the left
  and upper edge of the homotopy $I \times I \to Y \colon (s,t) \to
  h(u(s),t)$ into the lower and right edge
\begin{equation*}
  \xy 0;/r.15pc/: 
  (-20,-20)*{}="A"; (20,-20)*{}="B"; (-20,20)*{}="C"; (20,20)*{}="D";
   "A"; "B" **\dir{-}?(.5)*\dir{>}+(0,-4)*{f_0(u)} 
   \POS?(.25)*{}="-25" \POS?(.5)*{}="-50" \POS?(.75)*{}="-75";
   "C"; "D" **\dir{-}?(.5)*\dir{>}+(0,4)*{f_1(u)}
   \POS?(.25)*{}="+25" \POS?(.5)*{}="+50" \POS?(.75)*{}="+75";
    "A"; "C" **\dir{-}?(.5)*\dir{>}+(-8,0)*{h(x_0)}
    \POS?(.25)*{}="025" \POS?(.5)*{}="050" \POS?(.75)*{}="075";
    "B"; "D" **\dir{-}?(.5)*\dir{>}+(8,0)*{h(x_1)}
    \POS?(.75)*{}="125" \POS?(.5)*{}="150" \POS?(.25)*{}="175";
    "025"; "-25" **\dir{.}?(0.5)*\dir{>}; 
    "050"; "-50" **\dir{.}?(0.5)*\dir{>};
    "075"; "-75" **\dir{.}?(0.5)*\dir{>};
    "C"; "B" **\dir{.}?(0.5)*\dir{>}; 
    "+25"; "175" **\dir{.}?(0.5)*\dir{>};
    "+50"; "150" **\dir{.}?(0.5)*\dir{>};
    "+75"; "125" **\dir{.}?(0.5)*\dir{>};
   \endxy
\end{equation*}
we obtain a path homotopy $h(x_0) \cdot f_1(u) \simeq f_0(u) \cdot
h(x_1)$.
\end{proof}

\begin{cor}
   In the situation of Lemma~\ref{lemma:unbasedmaps}, the diagram
  \begin{equation*}
    \xymatrix@R=15pt{
      & {\pi_1(Y,f_1(x_0))} \ar[dd]_{\cong}^{[h(x_0)] - [\overline{h(x_0)}]} \\
      {\pi_1(X,x_0)} \ar[ur]^{(f_1)_*}\ar[dr]_{(f_0)_*} \\
      & {\pi_1(Y,f_0(x_0))} } 
  \end{equation*}
  commutes.
\end{cor}
\begin{proof}
  For any {\em loop\/} $u$ based at $x_0$, $f_0(u)h(x_0)=h(x_0)f_1(u)$
  or $f_0(u)=h(x_0)f_1(u)\overline{h(x_0)}$.
\end{proof}

\begin{cor}\label{cor:pi1hoeq}
  \begin{enumerate}
  \item  If \func{f}{X}{Y} is a  homotopy equivalence (possibly
    unbased) then the induced homo\m\
  \func{f_*}{\pi_1(X,x_0)}{\pi_1(Y,f(x_0))} is a group iso\m .
\item If \func{f}{X}{Y} is nullhomotopic (possibly unbased) then
  \func{f_*}{\pi_1(X,x_0)}{\pi_1(Y,f(x_0))} is the trivial homo\m .
  \end{enumerate}
  
\end{cor}
\begin{proof}
  (1) Let $g$ be a homotopy inverse to $f$ so that $gf \simeq 1_X$ and $fg
  \simeq 1_Y$. By Lemma~\ref{lemma:unbasedmaps} there is a commutative
  diagram
  \begin{equation*}
    \xymatrix{
      &&& {\pi_1(Y,f(x_0))} \\
      {\pi_1(X,x_0)} \ar@{=}[drr] \ar[r]^{f_*} 
      & {\pi_1(Y,f(x_0))} \ar@{=}[urr] \ar[r]_{g_*} 
      & {\pi_1(X,gf(x_0))} \ar[d]_{\cong} \ar[r]^{f_*} 
      & {\pi_1(Y,fgf(x_0))} \ar[u]^{\cong}\\
      && {\pi_1(X,x_0)} }
  \end{equation*}
  which shows that $g_*$ is both injective and surjective, ie $g_*$ is
  bijective. Then also $f_*$ is bijective.

  \noindent (2) If $f$ homotopic to a constant map $c$ then $f_*$
  followed by an iso\m\ equals $c_*$ which is trivial. Thuse also
  $f_*$ is trivial.
\end{proof}

We can now answer a question from Example~\ref{exmp:pi1Rn} and say
that any contractible space has trivial fundamental group.

\begin{defn}\label{prop:pathsin1con}
  A space is {\em simply connected\/} if  there is a unique path
  homotopy class between any two of its points.
\end{defn}

The space $X$ is simply connected if $\pi(X)(x_1,x_2)=\ast$ for all
$x_1,x_2 \in X$, or, equivalently, $X$ is path connected and
$\pi_1(X,x)=\ast$ at all points or at one point of $X$.

\section{Covering spaces}
\label{sec:covspaces}

A covering map over $X$ is a map that locally looks like the
projection map $X \times F \to X$ for some discrete space $F$. 

\begin{defn}
  A covering map is a \co\ surjective map \func{p}{Y}{X} with the
  property that for any point $x \in X$ there is a \n\ $U$ (an {\em
    evenly covered\/} \n ), a discrete set $F$, and a homeo\m\ $ U
  \times F \to p^{-1}(U)$ such that the diagram
  \begin{equation*}
    \xymatrix{
      U \times F \ar[dr]_{\mathrm{pr}_1} \ar[rr]^{\simeq} && 
      p^{-1}(U) \ar[dl]^{p \vert  p^{-1}(U)} \\
      & U }
  \end{equation*}
  commutes.
\end{defn}


Some covering spaces, but not all \eqref{sec:normal}, arise from left
group actions.  Consider a {\em left\/} action $G \times Y \to Y$ of a
group $G$ on a space $Y$.  Let \func{p_G}{Y}{G \backslash Y} be the
quotient map of $Y$ onto the orbit space $G \backslash Y$. The
quotient map $p_G$ is open because open subsets $U \subset Y$ have
open saturations $GU=\bigcup_{g \in G} gU = p_G^{-1}p_G(U)$
(\href{http://www.math.ku.dk/~moller/e03/3gt/notes/gtnotes.pdf}{General
  Topology} 2.82).  
The open sets in $\GY$ correspond bijectively to
saturated open sets in $Y$.

We now single out the left actions $G \times Y \to Y$ for which the
quotient map \func{p_G}{Y}{\GY} of $Y$ onto its orbit space is a
covering map.

\begin{defn}\label{defn:covspaceaction} \cite[(*) p.\ 72]{hatcher}
  A covering space action is a group action $G \times Y \to Y$ where
  any point $y \in Y$ has a \n\ $U$ such that the translated \n s
  $gU$, $g \in G$, are disjoint. (In other words, the action map $G
  \times U \to GU$ is a homeo\m .)
\end{defn}


\begin{exmp}
  The actions 
  \begin{itemize}
  \item $\Z\times \R \to \R \colon (n,t) \mapsto n+t$
  \item $\Z/2 \times S^n \to S^n \colon (\pm 1,x) \mapsto \pm x$
  \item $\Z/m \times S^{2n+1} \to S^{2n+1} \colon (\zeta,x) \mapsto
    \zeta x$, where $\zeta \in \C$ is an $m$th root of unity,
    $\zeta^m=1$,
  \item $\{\pm 1, \pm i, \pm j, \pm k \} \times S³ \to S³$, quaternion
    multiplication \cite[Example 1.43]{hatcher},
  \end{itemize}
  are covering space actions and the orbit spaces are $\Z \backslash
  \R =S¹$ (the circle), $\Z/2 \backslash S^n = \R P^n$ (real
  projective space), and $\Z/m \backslash S^{2n+1} = L^{2n+1}(m)$
  (lense space).  The action $\Z \times S¹\to S¹ \colon (n,z) \mapsto
  e^{\pi i \sqrt{2} n}z$ 
  is {\em not\/}
  a covering space action for the orbits are dense.
\end{exmp}

\begin{exmp}\label{exmp:cov}
  The maps \func{p_n}{S^1}{S^1}, $n \in \Z$, and
  \func{p_{\infty}}{\R}{S^1} given by $p_n(z)=z^n$, and
  $p_{\infty}(s)=e^{2\pi s} = (\cos(2\pi s), \sin(2\pi s))$ are
  covering maps of the circle with fibre $p_n^{-1}(1)=\Z/n\Z$ and
  $p_{\infty}^{-1}(1)=\Z$.  There are many covering maps of $S^1 \vee
  S^1$. The map $S^n \to C_2 \backslash S^n = \R P^n$, $n \geq 1$, is
  a covering map of real projective $n$-space.  The map $S^{2n+1} \to
  C_m \backslash S^{2n+1} = L^{2n+1}(m)$ is a covering map of the lens
  space. $M_g \to N_{g+1}$ a double covering map of the unorientable
  surface of genus $g+1$ with $F=\Z/2\Z$. Can you find a covering map
  of $M_g$? Can you find a covering map of $\R$?
\end{exmp}

\begin{thm}[Unique HLP for covering maps]\cite[1.30]{hatcher}\label{thm:HLP}
  Let \func{p}{Y}{X} be a covering map, $B$ be any space, and
  \func{h}{B \times I}{X} a homotopy into the base space.  If one end
  of the homotopy lifts to a map $B \times \{0\} \to Y$ then the whole
  homotopy admits a unique lift $B \times I \to Y$ such that the
  diagram
\begin{equation*}
  \xymatrix{
   B \times \{0\} \ar@{^(->}[d] \ar[r]^-{\widetilde{h}_0} & Y \ar[d]^p \\
   B \times I \ar[r]_h \ar@{-->}[ur]^{\widetilde{h}} & X }
\end{equation*}
commutes.
\end{thm}
\begin{proof}
  We consider first the case where $B$ is a point. The statement is
  then that  in the situation
  \begin{equation*}
    \xymatrix{
      {\{0\}} \ar[r]^{y_0} \ar@{^(->}[d] & Y \ar[d]^p \\
      I \ar[r]_u \ar@{-->}[ur]^{\widetilde{u}} & X}
  \end{equation*}
  there is a unique map $\widetilde{u} \colon I \to Y$ such that
  $p\widetilde{u}=u$ and $\widetilde{u}(0)=y_0$. For uniqueness of
  lifts from $I$ see Theorem~\ref{thm:lift}.(\ref{thm:lift1}). We need
  to prove existence.  The Lebesgue lemma
  (\href{http://www.math.ku.dk/~moller/e03/3gt/notes/gtnotes.pdf}{General
    Topology, 2.158}) 
  applied to the compact space $I$ says that there
  is a subdivision $0=t_0 < t_1 < \cdots < t_n=1$ of $I$ such that $u$
  maps each of the closed subintervals $[t_{i-1},t_i]$ into an evenly
  covered \n\ in $X$.  Suppose that we have lifted $u$ to
  $\widetilde{u}$ defined on $[0,t_{i-1}]$. Let $U$ be an evenly
  covered \n\ of $u(t_{i-1})$.  Suppose that the lift
  $\widetilde{u}(t_{i-1})$ belongs to $U \times \{\ell\}$ for some
  $\ell \in F$. Continue the given $\widetilde{u}$ with $(p \vert (U
  \times \{\ell\}))^{-1} \circ u \vert [t_{i-1},t_i]$.  After finitely
  many steps we have the unique lift on $I$.
  
  We now turn to the general situation. Uniqueness is clear for we
  have just seen that lifts are uniquely determined on the vertical
  slices $\{b\} \times I \subset B \times I$ for any point $b$ of $B$.
  Existence is also clear except that continuity of the lift is not
  clear.
  
  We now prove that the lift is \co . Let $b$ be any point of $X$. By
  compactness, there is a \n\ $N_b$ of $b$ and a subdivision $ 0=t_0 <
  t_1 < \cdots < t_n=1$ of $I$ such that $h$ maps each of the sets
  $N_b \times [t_{i-1},t_i]$ into an evenly covered \n\ of $X$.
  Suppose that $h(N_b \times [0,t_1])$ is contained in the evenly
  covered \n\ $U \subset X$ and let $\widetilde{U} \subset p^{-1}(U)
  \subset Y$ be a \n\ such that $p \vert \widetilde{U} \colon
  \widetilde{U} \to U$ is a homeo\m\ and $\widetilde{h}_0(b,0) \in
  \widetilde{U}$. We can not be sure that $\widetilde{h}_0(N_b \times
  \{0\}) \subset \widetilde{U}$; only if $N_b$ is connected.  Replace
  $N_b$ by $N_b \cap \widetilde{h}_0^{-1}(\widetilde{U})$. Then $
  \widetilde{h}_0(N_b \times \{0\}) \subset \widetilde{U}$. Then
  $(p\vert \widetilde{U})^{-1} \circ h\vert N_b \times [0,t_1]$ is a
  lift of $h \vert N_b \times [0,t_1]$ extending $\widetilde{h}_0$.
  After finitely many steps we have a lift defined on $N_b \times I$
  (where $N_b$ is possibly smaller than the $N_b$ we started with). Do
  this for every point $b$ of $B$.  These maps must agree on their
  overlap by uniqueness. So they define a lift $B \times I \to Y$.
  This lift is \co\ since it is \co\ on each of the open tubes $N_b
  \times I$.
\end{proof}

We emphasize the special case where $B$ is a point. Let $y_0 \in Y$ be
a point in $Y$ and $x_0=p(y_0) \in X$ its image in $X$.

\begin{cor}[Unique path lifting]
  Let $x_0$ and $x_1$ be two points in $X$ and let $y_0$ be a point in
  the fibre $p^{-1}(x_0) \subset Y$ over $x_0$. For any path \func uIX
  from $x_0$ to $x_1$, the exists a unique path \func{\tilde u}IY in
  $Y$ starting at $\tilde u(0)=y_0$. Moreover, homotopic paths have
  homotopic lifts: If \func vIX is a path in $X$
  that is path homotopic to $u$ then the lifts $\tilde u$ and $\tilde
  v$ are also path homotopic.
\end{cor}
\begin{proof}
  First, in Theorem~\ref{thm:HLP}, take $B$ to be point. Next, take
  $B$ to be $I$ and use the HLP to see that homotopic paths have
  homotopic lifts.
\end{proof}

\begin{cor} \label{eq:pi1actionsfibre} Let \func pYX be a covering map
  and let $y_0,y_1,y_2 \in Y$, $x_0=py_0, x_1=py_1, x_2=py_2$. 
  \begin{enumerate}
  \item By recording end points of lifts we obtain maps
\begin{equation*}
  p^{-1}(x_1) \times \pi(X)(x_1,x_2) \xrightarrow{\cdot}  p^{-1}(x_2), \qquad
  p^{-1}(x_0) \times \pi_1(X,x_0) \xrightarrow{\cdot} p^{-1}(x_0)
\end{equation*}
given by $y \cdot [u] = \widetilde{u}_y(1)$ where $\widetilde{u}_y$ is
the lift of $u$ starting at $y$.  
Multiplication by a path $u$ from $x_1$ to $x_2$ slides the fibre over
$x_1$ bijectively into the fibre over $x_2$.
\item The covering map \func pXY induces injective maps
 \begin{equation*}
    \pi(Y)(y_1,y_2) \xrightarrow{p_*} \pi(X)(x_1,x_2), \qquad 
    \pi_1(Y,y_0) \xrightarrow{p_*} \pi_1(X,x_0) 
  \end{equation*}
  The subset $p_*\pi(Y)(y_1,y_2) \subset \pi(X)(x_1,x_2)$ consists of
  all paths from $x_1$ to $x_2$ that lift to paths from $y_1$ to
  $y_2$. The subbgroup $p_*\pi_1(Y,y_0) \leq \pi_1(X,x_0)$ consists of
  all loops at $x_0$ that lloft to loops at $y_0$.
  \end{enumerate}
\end{cor}

\begin{defn}\label{eq:defnp-1}
The {\em monodromy functor\/} of the covering map \func pXY is a functor
\begin{equation*} F(p) \colon \pi(X) \to
  \mathbf{Set}
\end{equation*}
of the fundamental groupoid of the base space into the category
$\mathbf{Set}$ of sets. This functor takes a point in $x \in X$ to the
fibre $F(p)(x)=p^{-1}(x)$ over that point and it takes a path homotopy
class $u \in \pi(X)(x_0,x_1)$ to $F(p)(x_0)=p^{-1}(x_0) \to
p^{-1}(x_1)=F(p)(x_1) \colon y \to y \cdot u$. (The notation here is
such that $F(p)(uv) = F(p)(v) \circ F(p)(u)$ for paths $u \in
\pi(X)(x_0,x_1)$, $v \in \pi(X)(x_1,x_2)$.)
\end{defn}

In particular, the fibre $F(p)(x)=p^{-1}(x)$ over any point $x \in X$
is a right $\pi_1(X,x)$-set.

\begin{cor}[The fundamental groupoid of a covering space]
\label{cor:pi1Yinjects}
The fundamental groupoid of $Y$, 
\begin{equation*}
\pi(Y) = \pi(X) \rtimes F(p)  
\end{equation*}
is the Grothendieck construction of the fiber functor
\eqref{eq:defnp-1}.  In other words, the map $\pi(p) \colon
\pi(Y)(y_0,y_1) \to \pi(X)(x_0,x_1)$ is injective and the image is the
set of path homotopy classes from $x_0$ to $x_1$ that take $y_0$ to
$y_1$. In particular, the homo\m\
\func{p_*}{\pi_1(Y,y_0)}{\pi_1(X,x_0)} is injective and its image is
the set of loops at $x_0$ that lift to loops at $y_0$.
\end{cor}
\begin{proof}
  We consider the functor $F(p)$ as taking values in discrete categories.
  The objects of $\pi(X) \rtimes F(p)$ are pairs $(x,y)$ where $x \in X$
  and $y \in F(p)(x) \subset Y$. A \m\ $(x_1,y_1) \to (x_2,y_2)$ is a
  pair $(u,v)$ where $u$ is a \m\ in $\pi(X)$ from $x_1$ to $x_2$ and
  $v$ is a \m\ in $F(p)(x_2)$ from $F(p)(u)(x_1)= x_1 \cdot u$ to $y_2$. 
  As $F(p)(x_2)$ have no \m s but identities, the set of \m s $(x_1,y_1)
  \to (x_2,y_2)$ is the set of $u \in \pi(X)(x_1,x_2)$ such that $y_1
  \cdot u = y_2$. This is precisely $\pi(Y)(y_1,y_2)$.
\end{proof}

\begin{defn}
  For a space $X$, let $\pi_0(X)$ be the set of path components of $X$.
\end{defn}

\begin{lemma}\label{lemma:isotropy}
  Let \func pXY be a covering map.
  \begin{enumerate}
  \item Suppose that $X$ is path connected. The inclusion $p^{-1}(x_0)
    \subset Y$ induces a bijection $p^{-1}(x_0)/\pi_1(X,x_0) \to
    \pi_0(Y)$. In particular,
    \begin{equation*}
   \text{$Y$ is path connected} \iff 
    \text{$\pi_1(X,x_0)$ acts transitively on the fibre $p^{-1}(x_0)$}   
    \end{equation*}
    \item Suppose that $X$ and $Y$ are path connected. The maps
      \begin{alignat*}{2}
        \pi_1(Y,y_1) \backslash \pi(X)(x_1,x_2) &\leftrightarrow
        p^{-1}(x_2) &
        \qquad \pi_1(Y,y_0) \backslash \pi_1(X,x_0) &\leftrightarrow
        p^{-1}(x_0) \\
        \pi_1(Y,y_1)u &\to y_1 \cdot u &
           \pi_1(Y,y_0)u &\to y_0 \cdot u \\
           [pu_y] &\leftarrow y &
           [pu_y] &\leftarrow y
      \end{alignat*}
      are bijections. Here, $u_y$ is any path in $Y$ from $y_1$ or
      $y_0$ to $y$. In particular, $|\pi_1(X,x_0) \colon \pi_1(Y,y_0)|
      = |p^{-1}(x_0)|$.
  \end{enumerate}
 \end{lemma}
\begin{proof}
  The map $p^{-1}(x_0) \to \pi_0(Y)$, induced by the inclusion of the
  fibre into the total space, is onto because $X$ is path connected so
  that any point in the total space is connected by a path to a point
  in the fibre. Two points in the fibre are in the same path component
  of $Y$ if and only if are in the same $\pi_1(X,x_0)$-orbit.

  If $Y$ is path connected, then $\pi_1(X,x_0)$ acts transitively on
  the fibre $p^{-1}(x_0)$ with isotropy subgroup $\pi_1(Y,y_0)$ at
  $y_0$. 
\end{proof}


\begin{thm}[Lifting Theorem]\label{thm:lift}
  Let \func{p}{Y}{X} be a covering map and \func{f}{B}{X} a map into
  the base space. Choose base points such that $f(b_0)=x_0=p(y_0)$ and
  consider the lifting problem
  \begin{equation*}
    \xymatrix{
      & (Y,y_0) \ar[d]^p \\
      (B,b_0) \ar[r]_f \ar@{.>}[ur]^{\tilde f} & (X,x_0)}
  \end{equation*}
  \begin{enumerate}
  \item \label{thm:lift1} If $B$ is connected, then there exists at
    most one lift \func{\tilde{f}}{(B,b_0)}{(Y,y_0)} of $f$ over $p$.
  \item If $B$ is path connected and locally path connected then
    \begin{equation*}
      \text{There is a map \func{\tilde{f}}{(B,b_0)}{(Y,y_0)} such
        that $f = p\tilde f$} \iff  
    f_*\pi_1(B,b_0) \subset p_*\pi_1(Y,y_0) 
    \end{equation*}
      \end{enumerate}
\end{thm}
\begin{proof}
  (\ref{thm:lift1}) Suppose that $\widetilde{f}_1$ and
  $\widetilde{f}_2$ are lifts of the same map \func{f}{B}{X}. We claim
  that the sets $\{b \in B \mid \widetilde{f}_1(b) =
  \widetilde{f}_2(b)\}$ and $\{b \in B \mid \widetilde{f}_1(b) \neq
  \widetilde{f}_2(b)\}$ are open.

  Let $b$ be any point of $B$ where the two lifts agree.  Let $U
  \subset X$ be an evenly \n\ of $f(b)$. Choose $\widetilde{U} \subset
  p^{-1}(U) = U \times F$ so that the restriction of $p$ to
  $\widetilde{U}$ is a homoe\m\ and $\widetilde{f}_1(b) =
  \widetilde{f}_2(b)$ belongs to $\widetilde{U}$. Then
  $\widetilde{f}_1$ and $\widetilde{f}_2$ agree on the \n\
  $\widetilde{f}_1^{-1}(\widetilde{U}) \cap
  \widetilde{f}_2^{-1}(\widetilde{U})$ of $b$.

  Let $b$ be any point of $B$ where the two lifts do not agree.  Let
  $U \subset X$ be an evenly \n\ of $f(b)$. Choose disjoint open sets
  $\widetilde{U}_1,\widetilde{U}_2 \subset p^{-1}(U) = U \times F$ so
  that the restrictions of $p$ to $\widetilde{U}_1$ and
  $\widetilde{U}_2$ are homoe\m s and $\widetilde{f}_1(b)$ belongs to
  $\widetilde{U}_1$ and $\widetilde{f}_2(b)$ to $\widetilde{U}_2$.
  Then $\widetilde{f}_1$ and $\widetilde{f}_2$ do not agree on the \n\
  $\widetilde{f}_1^{-1}(\widetilde{U}_1) \cap
  \widetilde{f}_2^{-1}(\widetilde{U}_2)$ of $b$.

  \noindent (2) It is clear that if the lift exists, then the
  condition is satisfied. Conversely, suppose that the condition
  holds.  For any point $b$ in $B$, define a lift $\tilde f$ by
  \begin{equation*}
    \widetilde{f}(b)=y_0 \cdot [fu_b]
  \end{equation*}
  where $u_b$ is any path from $b_0$ to $b$.  (Here we use that $B$ is
  path connected.)  If $v_b$ is any other path from $b_0$ to $b$ then
  $y_0 \cdot [fu_b] = y_0 \cdot [fv_b]$ because $y_0 \cdot [fu_b \cdot
  f\overline{v_b}] = y_0$ as the loop $[fu_b \cdot f\overline{v_b}]
  \in \pi_1(Y,y_0)$ fixes the point $y_0$ by
  Lemma~\ref{lemma:isotropy}.

  We need to see that $\tilde{f}$ is \co . Note that any point $b \in
  B$ has a path connected \n\ that is mapped into an evenly covered
  \n\ of $f(b)$ in $X$. It is evident what $\tilde{f}$ does on this
  \n\ of $b$.
\end{proof}

A map \func{f}{B}{S^1 \subset \C-\{0\}} into the circle has an $n$th
root if and only if the induced homo\m\ \func{f_*}{\pi_1(B)}{\Z} is
divisible by $n$.

\section{The fundamental group of the circle, spheres, and lense spaces}
\label{sec:pi1S1}

For each $n \in \Z$, let $\omega_n$ be the loop
$\omega_n(s)=(\cos(2\pi ns), \sin(2\pi ns)$, $s \in I$, on the circle.

\begin{thm}\label{thm:pi1S1}
  The map $\func{\Phi}{\Z}{\pi_1(S^1,1)} \colon n \to [\omega_n]$
  is a group iso\m .
\end{thm}
\begin{proof}
  Let \func{p}{\R}{S^1} be the covering map $p(t)=(\cos(2\pi t),\sin(2
  \pi t))$, $t \in \R$. Remember that the total space $\R$ is simply
  connected as we saw in Example~\ref{exmp:pi1Rn}.  The fibre over $1$
  is $p^{-1}(1)=\Z$.  Let $u_n(t)=nt$ be the obvious
  path from $0$ to $n \in \Z$.  By Lemma~\ref{lemma:isotropy} the map
  \begin{equation*}
    \Z \to \pi_1(S^1,1) \colon n \to [pu_n] = [\omega_n]
  \end{equation*}
  is bijective.

  We need to verify that $\Phi$ is a group homo\m . Let $m$ and $n$ be
  integers. Then $u_m \cdot (m+u_n)$ is a path from $0$ to $m+n$ so
  it can be used instead of $u_{m+n}$ when computing $\Phi(m+n)$. We
  find that
  \begin{equation*}
   \Phi(m+n) = [p( u_m \cdot (m+u_n))] = [p( u_m) \cdot p(m+u_n)]
  = [p( u_m)][ p(m+u_n)] = [p( u_m)][ p(u_n)] = \Phi(m) \Phi(n) 
  \end{equation*}
     because $p(m+u_n) = pu_n$ as $p$ has period $1$. 
\end{proof}

\begin{thm}\label{thm:pi1Sn}
  The $n$-sphere $S^n$ is simply connected when $n>1$.
\end{thm}
\begin{proof}
  Let $N$ be the North and $S$ the South Pole (or any other two
  distinct points on $S^n$). The problem is that there are paths in
  $S^n$ that visit every point of $S^n$. But, in fact, any loop based
  at $N$ is {\em homotopic} to a loop that avoids $S$
  (\href{http://www.math.ku.dk/~moller/e02/3gt/opg/aug05.pdf}{Problem}
  and
  \href{http://www.math.ku.dk/~moller/e02/3gt/opg/ans.aug05.pdf}{Solution}).
  This means that $\pi_1(S^n-\{S\},N) \to \pi_1(S^n,N)$ is surjective.
  The result follows as $S^n-\{S\}$ is homeomorphic to the simply
  connected space $\R^n$.
\end{proof}

\begin{cor}
  The fundamental group of real projective $n$-space $\R P^n$ is
  $\pi_1(\R P^n) = C_2$ for $n>1$. The fundamental group of the lense space
  $L^{2n+1}(m)$ is $\pi_1(L^{2n+1}(m)) = C_m$ for $n>0$.
\end{cor}
\begin{proof}
  We proceed as in Theorem~\ref{thm:pi1S1}. Consider the case of the
  the covering map \func p{S^{2n+1}}{L^{2n+1}(m)} over the lense space
  $L^{2n+1}(m)$. Let $N=(1,0,\ldots,0) \in S^{2n+1} \subset \C^{n+1}$.
  The cyclic group $C_m = \gen{\zeta}$ of $m$th roots of unity is
  generated by $\zeta = e^{2\pi i/m}$. The map $\zeta^j \to \zeta^jN$,
  $j \in \Z$, is a bijection $C_m \to p^{-1}pN$ between the set $C_m$
  and the fibre over $pN$.  As $S^{2n+1}$ is simply connected there is
  a bijection
  \begin{equation*}
    \Phi \colon p^{-1}pN = C_m \to \pi_1(L^{2n+1}(m),pN) \colon \zeta^j \to
    [p\omega_j] 
  \end{equation*}
  where $\omega_j$ is the path in $S^{2n+1}$ from $N$ to $\zeta^jN$
  given by $\omega_j(s) = (e^{2\pi isj/m},0,\ldots,0)$. Since
  $\omega_{i+j} \simeq \omega_i \cdot (\zeta^i \omega_j)$, it follows
  just as in Theorem~\ref{thm:pi1S1} that $\Phi$ is a group homo\m .

  For the projective spaces, use the paths $\omega_j(s)=(\cos(2\pi js),
  \sin(2\pi j s),0,\ldots,0)$ from $N$ to $(-1)^jN$, to see that
 \begin{equation*}
    \Phi \colon p^{-1}pN = C_2 \to \pi_1(\R P^n,pN) \colon (\pm 1)^j \to
    [p\omega_j]
\end{equation*}
is a bijection. 
\end{proof}

\setcounter{subsection}{\value{thm}}
\subsection{Applications of $\pi_1(S^1)$}
\label{sec:applpi1S1}\add

Here are some standard applications of Theorem~\ref{thm:pi1S1}.

\begin{cor}\label{cor:degreenS1}
  The $n$th power homo\m\   $p_n \colon (S^1,1) \to (S^1,1) \colon z \to
  z^n$ induces the $n$th power homo\m\ $\pi_1(S^1,1) \to \pi_1(S^1,1)
  \colon [\omega] \to [\omega]^n$.
\end{cor}
\begin{proof}
  $(p_n)_*\Phi(1) = (p_n)_*[\omega_1] = [p_n \omega_1] = [\omega_1^n] =
  [\omega_n] =
  \Phi(n) = \Phi(1)^n$.
\end{proof}

\begin{thm}[Brouwer's fixed point theorem]\label{thm:brouwer}
  \begin{enumerate}
  \item The circle $S^1$ is not a retract of the disc $D^2$.
  \item Any map self-map of the disc $D^2$ has a fixed point.
  \end{enumerate}
\end{thm}
\begin{proof}
  (1)  Let
  \func{i}{S^1}{D^2} be the inclusion map. The induced map
  \func{i_*}{\Z=\pi_1(S^1)}{\pi_1(D^2)=0} is not injective so $S^1$
  can not be a retract by \ref{cor:piretract}. 

  \noindent (2) With the help of a fixed-point free self map of $D^2$
  one can construct a retraction of $D^2$ onto $S^1$. But they don't
  exist. 
\end{proof}


\begin{thm}[The fundamental theorem of algebra]
  Let  $p(z)=z^n + a_{n-1}z^{n-1} + \cdots + a_1z+a_0$ be a normed
  complex polynomial of degree $n$. If $n>0$, then $p$ has a root.
\end{thm}
\begin{proof}
  Any normed polynomial $p(z) = z^n + a_{n-1}z^{n-1} + \cdots +
  a_1z+a_0$ is nonzero when $|z|$ is large: When $|z| > 1+ |a_{n-1}| +
  \cdots + | a_0|$, then $p(z) \neq 0$ because
    \begin{multline*}
  \vert a_{n-1}z^{n-1} + \cdots + a_0 \vert 
   \leq  \vert a_{n-1} \vert \vert z \vert^{n-1} + \cdots + \vert a_0
  \vert 
  <\vert a_{n-1} \vert \vert z \vert^{n-1} + \cdots + \vert a_0
  \vert \vert z \vert^{n-1} \\ 
   = ( \vert a_{n-1} \vert  + \cdots + \vert a_0 \vert)
  \vert z \vert^{n-1} < \vert z \vert^n
  \end{multline*}
  Therefore any normed polynomial $p(z)$ defines a map $S^1(R) \to
  \C-\{0\}$ where $S^1(R)$ is the circle of radius $R$ and $R> 1+
  |a_{n-1}| + \cdots + | a_0|$. In fact, all the normed polynomials
  $p_t(z)=z^n + t(a_{n-1}z^{n-1} + \cdots + a_1z+a_0)$, $t \in I$,
  take $S^1(R)$ into $\C - \{0\}$ so that we have a homotopy
  \begin{equation*}
    S^1(R) \times I \to \C-\{0\} \colon (z,t) \to 
    z^n + t(a_{n-1}z^{n-1} + \cdots + a_1z+a_0)
  \end{equation*}
  between $p_1(z)=p(z) \vert S^1(R)$ and $p_0(z)=z^n$.

  If $p(z)$ has no roots at all, the map $p \vert S^1(R)$ factors
  through the complex plane $\C$ and is therefore nullhomotopic (as
  $\C$ is contractible) and so is the homotopic map $S^1(R) \to \C -
  \{0\} \colon z \to z^n$ and the composite map
  \begin{equation*}
    S^1 \xrightarrow{z \to Rz} 
    S^1(R) \xrightarrow{z \to z^n}
    \C - \{0\}  \xrightarrow{z \to z/|z|}
    S^1
  \end{equation*}
  But this is simply the map $S^1 \to S^1 \colon z \to z^n$ which we
  know induces multiplication by $n$ (\ref{cor:degreenS1}). However, a
  nullhomotopic map induces multiplication by $0$ (\ref{cor:pi1hoeq}).
  So $n=0$.
\end{proof}

A map \func{f}{S^1}{S^1} is {\em odd\/} if $f(-x)=-f(x)$ for all $x
\in S^1$. Any rotation (or reflection) of the circle is odd (because
it is linear).

\begin{lemma}
  Let \func{f}{S^1}{S^1} be an odd map. Compose $f$ with a rotation
  $R$ so that $Rf(1)=1$. The induced map
  \func{(Rf)_*}{\pi_1(S^1,1)}{\pi_1(S^1,1)} is multiplication by an
  odd integer. In particular, $f$ is not nullhomotopic.
\end{lemma}
\begin{proof}
  We must compute $(Rf)_*[\omega_1]$. 
  The HLP gives a lift
  \begin{equation*}
    \xymatrix{
      {\{0\}} \ar@{^(->}[d] \ar[rr]^0 && {\R} \ar[d] \\
      I \ar[r]_{\omega_1} \ar[urr]^{\widetilde{\omega}} 
      & S^1 \ar[r]_{Rf} & S^1}
  \end{equation*}
  and we have  $(Rf)_*[\omega_1]=[p\widetilde{\omega}]$. When $0 \leq
  s \leq 1/2$, $\omega_1(s+1/2) = -\omega_1(s)$ and also
  $Rf\omega_1(s+1/2) = -Rf\omega_1(s)$ as $Rf$ is odd. The lift,
  $\widetilde{\omega}$ of $Rf\omega_1$, then satisfies the equation
  \begin{equation*}
    \widetilde{\omega}(s+1/2) = \widetilde{\omega}(s) + q/2
  \end{equation*}
  for some {\em odd\/} integer $q$.  By continuity and connectedness
  of the interval $[0,1/2]$, $q$ does not depend on $s$. Now
  $\widetilde{\omega}(1) = \widetilde{\omega}(1/2)+q/2 =
  \widetilde{\omega}(0)+q/2+q/2=q$ and therefore $(Rf)_*[\omega_1] =
  [p\widetilde{\omega}]=[\omega_q]=[\omega_1]^q$. We conclude that
  $(Rf)_*$ is multiplication by the odd integer $q$. Since a
  nullhomopotic map induces the trivial group homomorphism
  \eqref{cor:pi1hoeq}, $f$ is not nullhomotopic.
\end{proof}

\begin{thm}[Borsuk--Ulam theorem for $n=2$]
  Let \func{f}{S^2}{\R^2} be any \co\ map. Then there exists a point
  $x \in S^2$ such that $f(x)=f(-x)$.
\end{thm}
\begin{proof}
  Suppose that  \func{f}{S^2}{\R^2} is a map such that $f(x) \neq
  f(-x)$ for all $x \in S^2$. The composite map
  \begin{equation*}
    \xymatrix@1@C=65pt{
     S^1 \ar@{^(->}[r]^{\textmd{incl}} & 
     S^2 \ar[r]^{x \to \frac{f(x)-f(-x)}{\vert f(x)-f(-x) \vert}} &
     S^1 }
  \end{equation*}
  is odd so it is not nullhomotopic. But the first map $S^1
  \hookrightarrow S^2$ is nullhomotopic because it factors through the
  contractible space $D^2_+=\{(x_1,x_2,x_3) \in S^2 \mid x_3 \geq
  0\}$. This is a contradiction.
\end{proof}

This implies that ther are no injective maps of $S^2 \to \R^2$; in
particular $S^2$ does not embed in $\R^2$.

\begin{prop}[Borsuk--Ulam theorem for $n=1$]
  Let \func{f}{S^1}{\R} be any \co\ map. Then there exists a point $x
  \in S^1$ such that $f(x)=f(-x)$.
\end{prop}
\begin{proof}
  Look at the map $g(x)=f(x)-f(-x)$. If $g$ is identically $0$,
  $f(x)=f(-x)$ for all $x \in S^1$. Otherwise, $g$ is an odd function,
  $g(-x)=-g(x)$, and $g$ has both positive and negative values. By
  connectedness, $g$ must assume the value $0$ at some point.
\end{proof}

This implies that there are no injective maps $S^1 \to \R$; in
particular $S^1$ does not embed in $\R$.

\section{The van Kampen theorem}
\label{sec:vankampen}

Let $G_j$, $j \in J$, be a set of groups indexed by the set $J$. The
{\em coproduct\/} (or {\em free product\/}) of these groups is a group
$\coprod_{i \in J}G_j$ with group homo\m s
\func{\varphi_j}{G_j}{\coprod_{j \in J} G_j} such that
\set\begin{equation}\label{eq:freeproduct}
  \Hom(\coprod_{j \in J} G_j,H) = \prod_{j \in J} \Hom(G_j,H) \colon
  \varphi  \to (\varphi \circ \varphi_j)_{j \in J}
\end{equation}\add
is a bijection for any group $H$.  The group $\coprod_{j \in J} G_j$
contains each group $G_j$ as a subgroup and these subgroups do not
commute with each other.  If the groups have presentations $G_j =
\gen{L_j \mid R_j}$ then $\coprod_{j \in J}\gen{L_j \mid R_j} =
\left\langle\cup_{j \in J} L_j \mid \cup_{j \in J}R_j \right\rangle$
as this group has the universal property.  See
\cite[6.2]{robinson:groups} for the construction of the free product.

The characteristic property \eqref{eq:freeproduct} applied to
$H=\prod_{j \in J}G_j$ shows that there is a group homo\m\ $\coprod
G_j \to \prod G_j$ from the free product to the direct product whose
restriction to each $G_j$ is the inclusion into the product.

\begin{exmp}\cite[Example II--III p $171$]{robinson:groups}
  \label{exmp:amalg} 
  $\Z/2 \amalg \Z/2 = \Z \rtimes \Z/2$ and $\Z/2 \amalg \Z/3 =
  \mathrm{PSL}(2,\Z)$. We can prove the first assertion:
  \begin{equation*}
    \Z/2 \amalg \Z/2 = \langle a,b \mid a^2,b^2 \rangle
    = \langle a,b,c \mid a^2,b^2,c=ab \rangle
    = \langle a,b \mid a^2,acac,c \rangle 
    = \langle a,b \mid a^2,c,aca=c^{-1} \rangle
  \end{equation*}
  but the second one is more difficult.
\end{exmp}

Suppose that the space $X= \bigcup_{j \in J}X_j$ is the union of open
and path connected subspaces $X_{j}$ and that $x_0$ is a point in
$\bigcap_{j \in J} X_{j}$. 
The inclusion of the subspace $X_j$ into $X$ induces a group homo\m\
\func{\iota_j}{\pi_1(X_j,x_0)}{\pi_1(X,x_0)}.
The coproduct $\coprod_{j \in J}
\pi_1(X_j,x_0)$ is a group equipped with group homo\m s
\func{\varphi_j}{\pi_1(X_j,x_0)}{\coprod_{j \in J} \pi_1(X_j,x_0)}. Let
  \begin{equation*}
    \func{\Phi}{\amalg_{j \in J} \pi_1(X_j,x_0)}{\pi_1( \bigcup_{j
        \in J}X_j,x_0) = \pi_1(X,x_0)}
  \end{equation*}
  be the group homo\m\ 
  determined by $\Phi  \circ \varphi_j = \iota_j$.

  Is $\Phi$ surjective? In general, no. The circle, for instanec, is
  the union of two contractible open subspaces, so $\Phi$ is not onto
  in that case. But, if any loop in $X$ is homotopic to a product of
  loops in one of the subspaces $X_j$, then $\Phi$ is surjective.

  Is $\Phi$ injective?  It will, in general, not be injective, because
  the individual groups $\pi_1(X_i)$ in the free product do not
  intersect but the subspaces do intersect.  Any loop in $X$ that is a
  loop in $X_i \cap X_j$ will in the free product count as a loop both
  in $\pi_1(X_i)$ and in $\pi_1(X_j)$. We always have commutative
  diagrams of the form
    \begin{equation*}
      \xymatrix{
        & {\pi_1(X_i,x_0)} \ar[dr]^{\iota_i} \\
       {\pi_1(X_i \cap X_j,x_0)} \ar[ur]^{\iota_{ij}}
       \ar[dr]_{\iota_{ji}} 
       && { \pi_1(X,x_0)} \\
       &  {\pi_1(X_j,x_0)} \ar[ur]_{\iota_j} }
    \end{equation*}
    where $\iota_{ij}$ are inclusion maps. This means that
    $\Phi(\iota_{ij}g)=\Phi(\iota_{ji}g)$ for any $g \in
    \pi_1(X_i \cap X_j,x_0)$ so that 
    \set\begin{equation}\label{eq:kerPhi}
      \forall i,j \in J \forall g \in \pi_1(X_i \cap X_j) \colon 
       {\iota_{ij}}(g) {\iota_{ji}}(g)^{-1} \in \ker \Phi
    \end{equation}\add
    Let $N \leq \coprod_{j \in J} \pi_1(X_j,x_0)$ be the smallest
    normal subgroup containing all the elements of \eqref{eq:kerPhi}.
    The kernel of $\Phi$ must contain $N$ but, of course, the kernel
    could be bigger.  The surprising fact is that often it isn't.

\begin{thm}[Van Kampen's theorem]\label{thm:vankampen}
  Suppose that $X= \bigcup_{j \in J}X_j$ is the union of open and
  path connected subspaces $X_{j}$ and that $x_0$ is a point in
  $\bigcap_{j \in J} X_{j}$.
  \begin{enumerate}
  \item If the intersection of any two of the open subspaces is
    path connected then $\Phi$ is surjective.
  \item If the intersection of any three of the open subspaces is
    path connected then the kernel of $\Phi$ is $N$.
  \end{enumerate} 
  \end{thm}

  \begin{cor}
    If the intersection of any three of the open subspaces is
    path connected
     then $\Phi$ determines an iso\m\
  \begin{equation*}
  {\overline{\Phi}} \colon 
  {\coprod_{j \in J} \pi_1(X_j,x_0)}/N \cong \pi_1(X,x_0)   
  \end{equation*}
  \end{cor}

\begin{proof}[Proof of Theorem~\ref{thm:vankampen}]
  (1) We need to show that any loop $u \in \pi_1(X)$ in $X$ is a
  product $u_1 \cdots u_m$ of loops $u_i \in \pi_1(X_{j_i})$ in one of
  the subspaces.  Let \func{u}{I}{X} be a loop in $X$. 
  
  Thanks to the Lebesgue lemma
  (\href{http://www.math.ku.dk/~moller/e03/3gt/notes/gtnotes.pdf}{General
    Topology, 2.158}) 
  we can find a subdivision $0=t_0 <t_1 < \cdots
  t_m = 1$ of the unit interval so that $u_i=u|[t_{i-1},t_i]$ is a
  path in (say) $X_i$. As $u(t_i) \in X_{i} \cap X_{i+1}$, and also
  the base point $x_0 \in X_{i} \cap X_{i+1}$, and $X_i \cap X_{i+1}$
  is path connected, there is path $g_i$ in $X_{i} \cap X_{i+1}$ from
  the basepoint $x_0$ to $u(t_{i-1})$. The situation looks like this:
  \begin{equation*}
  \xy 0;/r.15pc/: 
   (0,0)*{}="A"; (-20,20)*{}="B"; (20,20)*{}="C";
   (-30,30)*{X_1}; 
   (30,30)*{X_2}; 
   (-10,-10)*{X_3}; 
   (0,12)*{\ast}="O"; 
   (0,25)*{\bullet}="O12"; 
   (15,10)*{\bullet}="O23"; 
   "A" *\xycircle(30,20){.};
   "B" *\xycircle(30,20){.};
   "C" *\xycircle(30,20){.};
   "O"; "O12" **\dir{--}?(.5)*\dir{>}+(4,-2)*{g_1};;
   "O"; "O12" **\crv{(-40,-10) & (-40,40)}?(.55)*\dir{>}+(-12,0)*{u|[0,t_1]};
   "O12"; "O23" **\crv{(40,40) & (40,-10)}?(.5)*\dir{>}+(12,0)*{u|[t_1,t_2]};
   "O"; "O23" **\dir{--}?(.6)*\dir{>}+(2,-4)*{g_2};
   "O23"; "O" **\crv{(20,-30) & (-10,-10)}?(.5)*\dir{>}+(0,-5)*{u|[t_2,1]};
    \endxy
\end{equation*}
Now $u \simeq u|[0,t_1] \cdot u|[t_1,t_2] \cdots u|[t_{m-1},1] \simeq
(u|[0,t_1] \cdot \overline{g_1}) \cdot (g_1 \cdot u|[t_1,t_2] \cdot
\overline{g_2}) \cdots (g_m \cdot u|[t_{m-1},1])$ is a product of
loops where each factor is inside one of the subspaces.

\noindent (2) 
Let $N \triangleleft \amalg \pi_1(X_i)$ be the smallest normal subgroup
containing all the elements (\ref{eq:kerPhi}).  Let $u_i \in
\pi_1(X_{j_i})$. For simplicity, let's call $X_{j_i}$ for $X_i$.
Consider the product
\begin{equation*}
  \underbrace{u_1}_{\pi_1(X_1)} \underbrace{u_2}_{\pi_1(X_2)} \cdots 
  \underbrace{u_m}_{\pi_1(X_m)} \in \amalg_{j \in J} \pi_1(X_j)
\end{equation*}
and suppose that $\Phi(u_1 \cdots u_m)$ is the unit element of
$\pi_1(X)$.  We want to show that $u_1 \cdots u_m$ lies in the normal
subgroup $N$ or that $u_1 \cdots u_m$ is the identity in the quotient
group $\amalg \pi_1(X_j)/N$.

Since $u_1 \cdots u_m$ is homotopic to the constant loop in $X$ there
is homotopy $I \times I \to X=\bigcup X_j$ from the loop $u_1 \cdots
u_m$ in $X$ to the constant loop. Divide the unit square $I \times I$
into smaller rectangles such that each rectangle is mapped into one of
the subspaces $X_j$. We may assume that the subdivision of $I \times
\{0\}$ is a further subdivision of the subdivision at $i/m$ coming
from the product $u_1 \cdots u_m$.  It could be that one new vertex is
(or more new vertices are) inserted between $(i-1)/m$ and $i/m$.
\begin{equation*} 
  \xy 
  (-20,0)*{}="A"; (20,0)*{}="B";
  "A"; "B" **\dir{-}
   \POS?(0)*{}+(0,-4)*{\ast} 
   \POS?(1)*{}+(0,-4)*{\ast} 
   \POS?(.4)*{\bullet}
   \POS?(.2)*{}+(0,4)*{X_k}
   \POS?(.6)*{}+(0,4)*{X_{\ell}};
   \POS?(.5)*{}+(0,-4)*{X_{i}};
  \endxy
\end{equation*} 
Connect the image of the new vertex $\bullet$ with a path $g$ inside
$X_{i} \cap X_k \cap X_{\ell}$ to the base point. Now $u_i$ is
homotopic in $X_i$ to the product $(u_i|[(i-1)/m,\bullet] \cdot
\overline{g})\cdot (g \cdot u_i|[\bullet,i/m])$ of two loops in $X_i$.
This means that we may as well assume that no new subdivision points
have been introduced at the bottom line $I \times \{0\}$.  Now perturb
slightly the small rectangles, but not the ones in the bottom and top
row, so that also the corner of each rectangle lies in at most three
rectangles. The lower left corner may look like this:
\begin{equation*} 
  \xy 0;/r.20pc/:  
  (-20,0)*{}="A"+(-4,0)*{\ast}; (20,0)*{}="B";  (20,-20)*{}="C"+(0,-4)*{\ast}; 
  (30,20)*{}="D";  
  (30,0)*{}="E"; (-20,20)*{}="F"; (-20,-20)*{}="G"+(-4,-4)*{\ast}; 
  (30,-20)*{}="H";
  "A"; "B" **\dir{-}
   \POS?(1)*{\triangle}="trekant"
   \POS?(.7)*{\bullet}="prik"
   \POS?(.2)*{}+(5,10)*{X_5}
   \POS?(.9)*{}+(5,10)*{X_6}
   \POS?(.5)*{}+(0,-10)*{X_1}
   \POS?(.35)*{}+(0,-4)*{u_{15}}
   \POS?(.85)*{}+(0,-4)*{u_{16}};
   (30,-10)*{X_2};
   "B"; "E" **\dir{-};
   "B"; "C" **\dir{-} \POS?(.5)*{}+(-4,0)*{u_{12}};
   (-20,20);(20,20) **\dir{}
    \POS?(1)*{}="utrekant"
    \POS?(.7)*{}="uprik";
    "prik"; "uprik" **\dir{-};
    "F"; "A" **\dir{-};
    "G"; "A" **\dir{-};
    "G"; "C" **\dir{-} \POS?(.5)*{}+(0,-4)*{u_1};
    "F"; "D" **\dir{-};
    "C"; "H" **\dir{-};
  \endxy
\end{equation*}
The loop $u_1$ in $X_1$ is homotopic to the product of paths
$u_{15}u_{16}u_{12}$ by a homotopy as in the proof of
\ref{lemma:unbasedmaps}. Connect the image of the point $\bullet$ to
the base point by a path $g_{156}$ inside $X_1 \cap X_5 \cap X_6$ and
connect the image of the point $\triangle$ to the base point by a path
$g_{126}$ inside $X_1 \cap X_2 \cap X_6$.  Then $u_1$ is homotopic in
$X_1$ to the product of loops $(u_{15}\overline{g_{156}}) \cdot
(g_{156}u_{16}\overline{g_{126}}) \cdot (g_{126}u_{12})$ in $X_1$. The
first of these loops is a loop in $X_1 \cap X_5$, the second is a loop
in $X_1 \cap X_6$, and the third is a loop in $X_1 \cap X_2$. In $\amalg
\pi_1(X_j)$ and  modulo
the normal subgroup $N$ we have that 
\begin{equation*}
  \underbrace{u_1}_{X_1} \underbrace{u_2}_{X_2} \cdots =
  \underbrace{u_{15}\overline{g_{156}}}_{X_1} \cdot
  \underbrace{g_{156}u_{16}\overline{g_{126}}}_{X_1} \cdot
  \underbrace{g_{126}u_{12}}_{X_1}\underbrace{u_2}_{X_2}\cdots =
  \underbrace{u_{15}\overline{g_{156}}}_{X_5} \cdot
  \underbrace{g_{156}u_{16}\overline{g_{126}}}_{X_6}\cdot
  \underbrace{g_{126}u_{12} \cdot u_2}_{X_2}\cdots
\end{equation*}
After finitely many steps we conclude that modulo $N$ the product $u_1
\cdots u_m$ equals a product of constant loops, the identity element.
\end{proof}

\begin{cor}\label{prop:pi1ofwedge}
  Let $X_j$ be a set of path connected spaces. Then
  \begin{equation*}
    \coprod_{j \in J} \pi_1(X_j) \cong \pi_1(\bigvee_{j \in J} X_j)
  \end{equation*}
  provided that each base point $x_j \in X_j$ is the deformation
  retract of an open \n\ $U_j \subset X_j$.
\end{cor}
\begin{proof}
  Van Kampen's theorem does not apply directly to the subspaces $X_j$
  of $\bigvee X_j$ because they are not open. Instead, let $X_j'= X_j
  \cup \bigvee_{i \in J} U_i$. The subspaces $X_j'$ are open and path
  connected and the intersection of at least two of them is the
  contractible space $\bigvee_{i \in J} U_i$.  Moreover, $X_j$ is a
  deformation retract of $X_j$.
\end{proof}

For instance, punctured compact surfaces have free fundamental groups.


\begin{cor}[van Kampen with two subspaces]\label{cor:amalg}
  Suppose that $X=X_1 \cup X_2$ where $X_1$, $X_2$, and $X_1 \cap X_2
  \neq \emptyset$ are open and path connected. Then
  \begin{equation*}
    \pi_1(X_1 \cup X_2,x_0) \cong
\pi_1(X_1,x_0) \amalg_{\pi_1(X_1 \cap X_2,x_0)} \pi_1(X_2,x_0) \cong
    \end{equation*}
  for any basepoint $x_0 \in X_1 \cap X_2$.
\end{cor}

This means that when $X_1 \cap X_2$ is path connected the fundamental
group functor takes a push out of spaces to a push out, amalgamated
product, of groups
\begin{equation*}
  \xymatrix@R=10pt{
     X_1 \cap X_2 \ar[r]^-{i_1} \ar[dd]_{i_2} &
     X_2 \ar[dd]  &&& {\pi_1(X_1 \cap X_2)}  \ar[r]^-{(i_1)_*}
     \ar[dd]_{(i_2)_*} & {\pi_1(X_2)} \ar[dd] \\ 
      && {} \ar[r]^{\pi_1} & {} \\
     X_2 \ar[r] & X &&& {\pi_1(X_2)} \ar[r] & {\pi_1(X)} }
\end{equation*}
As a very special case, we see that a space, that is the union of
two open simply connected subspaces with path connected intersection,
is simply connected. This proves, again (Theorem~\ref{thm:pi1Sn}),
that $S^n$ is simply connected when $n > 1$. 

We can use this simple variant of van Kampen to analyze the effect on
the fundamental group of attacing cells.

\begin{cor}[The fundamental group of a cellular
  extension]\label{prop:pi1plusncell}
  Let $X$ be a path connected space. Then
  \begin{equation*}
   \pi_1(X \cup_{\coprod f_{\alpha}} \coprod D^n_{\alpha}) = 
   \begin{cases}
     \pi_1(X)/\gen{\gamma_{\alpha} f_{\alpha} \overline{\gamma_{\alpha}}}
     & n=2 \\
    \pi_1(X) & n>2
   \end{cases}
  \end{equation*}
  where $\gamma_{\alpha}$ is a path from the base point of $X$ to the
  image of the base point of $S^{1}_{\alpha} \subset D^2_{\alpha}$.
\end{cor}
\begin{proof} 
  Let $Y$ be $X$ with the $n$-cells attached. Attach strips, fences
  connecting the base point of $X$ with the base points of the
  attached cells, to $Y$ and call the results $Z$. This does not
  change the fundamental group as $Y$ is a deformation retract of $Z$
  (Corollary~\ref{cor:piretract}). Let $A$ be $Z$ with the top half of
  each cell removed and let $B=Z-X$. Then $Z=A \cup B$ and $A \cap B$
  are path connected (the fences are there to make $A$ and $B$ path
  connected) so that
  \begin{equation*}
    \pi_1(Z) = \pi_1(A) \amalg_{\pi_1(A \cap B)} \pi_1(B)
  \end{equation*}
  by the van Kampen theorem in the simple form of
  Corollary~\ref{cor:amalg}. Now $B$ is contractible, hence simply
  connected (Corollary~\ref{cor:pi1hoeq}), so $\pi_1(Y)=\pi_1(Z)$ is
  the quotient of $\pi_1(A)$ by the smallest normal subgroup
  containing the image of $\pi_1(A \cap B) \to \pi_1(A)$. But $A \cap
  B$ is homotopy equivalent to a wedge
  $\bigvee_{\alpha}S^{n-1}_{\alpha}$ of $(n-1)$-spheres. In
  particular, $A \cap B$ is simply connected when $n>2$
  (Corollary~\ref{prop:pi1ofwedge}, Theorem~\ref{thm:pi1Sn}) so that
  $\pi_1(Y) = \pi_1(Z) = \pi_1(A) = \pi_1(X)$.  When $n=2$, $\pi_1(A
  \cap B)$ is a free group and the image of it in $\pi_1(A)=\pi_1(X)$
  is generated by the path homotopy classes of the loops
  $\gamma_{\alpha}f_{\alpha}\overline{\gamma_\alpha}$.
\end{proof}

\begin{cor}\label{cor:pi1CW}
  Let $X$ be a CW-complex with skeleta $X^k$, $k \geq 0$. Then
  \begin{equation*}
    \pi_0(X^1) = \pi_0(X), \qquad\qquad \pi_1(X^2)=\pi_1(X)
  \end{equation*}
\end{cor}

\begin{cor}
 The fundamental groups of the compact surfaces of positive genus $g$ are 
 \begin{equation*}
   \pi_1(M_g) = \langle a_1,b_1,\ldots ,a_g,b_g \mid \prod [a_i,b_i]
   \rangle, \qquad
    \pi_1(N_g) = \langle a_1,\ldots ,a_g \mid \prod a_i^2
   \rangle, \qquad
 \end{equation*}
\end{cor}

The compact orientable surfaces $M_g$, $g \geq 0$, are distinct,
$\pi_1(M_g)_{\mathrm{ab}}=\Z^{2g}$, and the compact nonorientable
surfaces $N_h$, $h \geq 1$, are distinct,
$\pi_1(N_g)_{\mathrm{ab}}=\Z^{g} \times \Z/2$.

\begin{cor}\label{cor:pi1plusncell}
  Let $M$ be a connected manifold of dimension $ \geq 3$. Then
  $\pi_1(M-\{x\})=\pi_1(M)$ for any point $x \in M$.
\end{cor}
\begin{proof}
  Apply van Kampen to $M = M-\{x\} \cup D^n$, $M-\{x\} \cap D^n \simeq
  S^{n-1}$ and remember that $S^{n-1}$ is simply connected when $n
  \geq 3$.
\end{proof}

Which groups can be realized as fundamental groups of spaces?  For
instance, $C_\infty=S^1$ and $C_m = S^1 \cup_m D^2$ so that any
finitely generated abelian group can be realized as the fundamental
group of a product of these spaces.

\begin{cor}\label{cor:XG}
  For any group $G$ there is a $2$-dimensional CW-complex $X_G$ such
  that $\pi_1(X_G) \cong G$.
\end{cor}
\begin{proof}
  Choose a presentation $G = \langle g_{\alpha} \mid
  r_{\beta} \rangle$ and let
\begin{equation*}
  X_G = D^0 \cup \coprod_{\{g_{\alpha}\}} D^1
            \cup \coprod_{\{r_{\beta}\}} D^2 
\end{equation*}
be the $2$-dimensional CW-complex whose $1$-skeleton is a wedge of
circles, one for each generator, with $2$-discs attached along the
relations.
\end{proof}

Observe that $X_{H \amalg G} = X_H \vee X_G$. Also,
$X_{\pi_1(M_g)}=M_g$, $X_{\pi_1(N_g)}=N_g$, $g \geq 1$.

\subsection{Fundamental groups of knot and link complements}
\label{sec:pi1knot}

The complement of a pair of unlinked circles in $\R^3$ deformation
retracts to $S^1 \vee S^1 \vee S^2 \vee S^2$ and a pair of linked
circles to $(S^1 \times S^1) \vee S^2$. The fundamental groups are $\Z
\ast \Z$ and $\Z \times \Z$, respectively. Thus the two complements
are not homeomorphic.

Let $m$ and $n$ be relatively prime natural numbers and $K=K_{mn}$ the
$(m,n)$-torus knot.  We want to compute the knot group
$\pi_1(\R^3-K)$. 

According to \eqref{cor:pi1plusncell}, $\pi_1(\R^3-K)=\pi_1(S^3-K)$. Now
\begin{equation*}
  S^3 = \partial D^4 = \partial (D^2 \times D^2) = \partial D^2 \times
  D^2 \cup D^2 \times \partial D^2
\end{equation*}
is the union of two solid tori  intersecting in a torus $S^1 \times
S^1$. Let $K$ be embedded in this middle torus. Then
\begin{equation*}
  S^3-K = ( \partial D^2 \times
  D^2 -K) \cup (D^2 \times \partial D^2 -K), \qquad
  ( \partial D^2 \times
  D^2 -K) \cap (D^2 \times \partial D^2 -K) = S^1 \times S^1 -K
\end{equation*}
and van Kampen says (if we ignore\footnote{To fix this, thicken the
  knot and enlarge the two solid tori a little so that they overlap.} 
the condition that the subsets
should be open)
\begin{equation*}
  \pi_1(S^3-K)= \frac{ \pi_1( \partial D^2 \times
  D^2 -K) \amalg \pi_1(D^2 \times \partial D^2 -K) }
  { \pi_1( S^1 \times S^1 -K)}
\end{equation*}
Here, $\partial D^2 \times D^2 -K$ deformation retracts onto the core
circle $\partial D^2 \times \{0\}$, and $S^1 \times S^1 -K$ (the torus
minus the knot) is an annulus $S^1 \times (0,1)$. (Take an open strip
$[0,1] \times (O,1)$ and wrap it around the torus so that the end $0
\times (0,1)$ meets the end $1 \times (0,1)$). The image of the
generator of this infinite cyclic group is the $m$ power of a
generator, respectively the $n$th power. Hence
  \begin{equation*}
     \pi_1(S^3-K) = \gen{a,b \mid a^m=b^n} = G_{mn}
  \end{equation*}

  It is now a matter of group theory to tell us that if $G_{m_1n_1}$
  and $G_{m_2n_2}$ are isomorphic then $\{m_1,n_1\} = \{m_2,n_2\}$.
  In order to analyze this group, note that $a^m=b^n$ is in the
  center. Let $C$ be the central group generated by this element. The
  quotient group
  \begin{equation*}
    G_{mn}/C =  \gen{a,b \mid a^m, b^n} = \Z/m \amalg \Z/n
  \end{equation*}
  has no center. (In general the free product $G\amalg H$ of two
  nontrivial groups has no center because the elements are words in
  elements from $G$ alternating with elements from $H$.) Therefore $C$
  is precisely the center of $G_{mn}$. Thus we can recover $mn$ as the
  order of the abelianization of $G/Z(G)$. Also, any element of finite
  order in $ \Z/m \amalg \Z/n$ is conjugate to an element of $\Z/m$ or
  $\Z/n$. Thus we can recover the largest of $m,n$ as the maximal
  order of a torsion element in $G/Z(G)$. Thus we can recover the set
  $\{m,n\}$. 
  \begin{cor}
    There are infinitely many knots. 
    (\href{http://www.pims.math.ca/knotplot/zoo/}{Here} are some of them.) 
  \end{cor}

  Another way of saying this is that $\partial D^2 \times
  D^2 -K$ deformation retracts onto the mapping cylinder of the degree
  $m$, respectively $n$, map $S^1 \to S^1$. Thus the union of these
  two spaces, $S^3-K$, deformation retracts onto the union of the two
  mapping cylinders, which is the double mapping cylinder $X_{mn}$ for
  the two maps.

  Thus $X_{mn}$ embeds in $S^3$ and $\R^3$ when $(m,n)=1$. On the other
  hand $X_{22}$ is the union of two M\"obius bands. A  M\"obius band
  is $\R P^2$ minus an open $2$-disc, so $X_{22} = \R P^2 \# \R P^2$,
  the Klein bottle, which does not embed in $\R^3$.

\section{Categories}
\label{sec:cat}

A {\em category\/}  $\mathcal C$ consists of \cite{maclane:cat}
\begin{itemize}
\item {\em Objects\/} $a,b,\ldots$
\item For each pair of objects $a$ and $b$ a set of {\em \m s\/} $\mathcal
  C(a,b)$ with domain $a$ and codomain $b$
\item A composition function $\mathcal C(b,c) \times \mathcal C(a,b)
  \to \mathcal C(a,c)$ that to each pair of \m s $g$ and $f$ with
  $\mathrm{dom}(g)=\mathrm{cod}(f)$ associates a \m\ $g \circ f$ with
  $\mathrm{dom}(g \circ f)= \mathrm{dom}(f)$ and $\mathrm{cod}(g \circ
  f)= \mathrm{cod}(g)$ 
\end{itemize}
We require 
\begin{description}
\item[Identity] 
For each object $a$ the \m\ set $\mathcal C(a,a)$
contains a \m\ $\mathrm{id}_a$ such that $g \circ \mathrm{id}_a = g$
and $\mathrm{id}_a \circ f = f$ whenever these compositions are
defined
\item[ Associativity ] $h \circ (g \circ f) = (h \circ g) \circ f$
  whenever these compositions are defined
\end{description}

A \m\ $f \in \mathcal C(a,b)$ with domain $a$ and codomain $b$ is
sometimes written \func{f}{a}{b}. A \m\ \func{f}{a}{b} is an {\em
  iso\m\/} if there exists a \m\ \func{g}{b}{a} such that the two
possible compositions are the respective identities.

\begin{defn}
  A group is a category with one object where all \m s are iso\m s. A
  groupoid is a category where all \m s are iso\m s.
\end{defn}
 
\begin{exmp}
  In the category $\mathbf{Top}$ of topological spaces, the objects
  are topological spaces, the \m s are \co\ maps, and composition is
  the usual composition of maps.  In the category $\mathbf{hoTop}$,
  the objects are topological spaces, the \m s are homotopy classes of
  \co\ maps, and composition is induced by the usual composition of
  maps. In the category $\mathbf{Grp}$ of groups, the objects are
  groups, the \m s are groups homo\m s, and composition is the usual
  composition of group homo\m s. In the category $\mathbf{Mat}_R$ the
  objects are the natural numbers $\Z_+$, the set of \m s $m \to n$
  consists of all $n$ by $m$ matrices with entries in the commutative
  ring $R$, and composition is matrix multiplication. The fundamental
  groupoid $\pi(X)$ of a topological space $X$ is a groupoid where the
  objects are the points of $X$ and the \m s $x \to y$ are the
  homotopy classes $\pi(X)(x,y)$ of paths from $x$ to $y$, and
  composition is composition of path homotopy classes.
\end{exmp}

A {\em functor\/} \func{F}{\mathcal C}{\mathcal D} associates to each
object $a$ of $\mathcal C$ an object $F(a)$ of $\mathcal D$ and to
each \m\ \func{f}{a}{b} in $\mathcal C$ a \m\ \func{F(f)}{F(a)}{F(b)}
in $\mathcal D$ such that $F(\mathrm{id}_a)=\mathrm{id}_{F(a)}$ and
$F(g \circ f)=F(g) \circ F(f)$.

A {\em natural transformation\/} $\tau \colon F \Longrightarrow G
\colon \mathcal{C} \to \mathcal{D}$ between two functors
\func{F,G}{\mathcal{C}}{\mathcal{D}} is a $\mathcal{D}$-\m\ $\tau(a) \in
\mathcal{D}(Fa,Ga)$ for each object $a$ of $\mathcal{C}$ such that the
diagrams
\begin{equation*}
  \xymatrix{
    a \ar[d]_f & Fa \ar[d]_{Ff} \ar[r]^{\tau(a)} & Ga
    \ar[d]^{Gf}\\
    b & Fb \ar[r]^{\tau(b)} & Gb}
\end{equation*}
commute for all \m s $f \in \mathcal{C}(a,b)$ in $\mathcal{C}$.  A
natural transformation $\tau$ is a {\em natural iso\m\/} if all the
components $\tau(a)$, $a \in \mathrm{Ob}(\mathcal{C})$, are
$\mathcal{D}$-iso\m s.

\begin{exmp}
  The fundamental group is a functor from the category of based
  topological spaces and based homotopy classes of maps to the
  category of groups. 

  The fundamental groupoid is a functor from the category of
  topological spaces to the category of groupoids.  Any homotopy $h
  \colon f_0 \simeq f_1$ induces a natural iso\m\ $h \colon \pi(f_0)
  \Longrightarrow \pi(f_1) \colon \pi(X) \to \pi(Y)$ between functors
  between fundamental groupoids (Lemma~\ref{lemma:unbasedmaps}).
\end{exmp}

\begin{defn}
  Let $\mathcal{C}$ and $\mathcal{D}$ be categories. The functor
  category $\mathrm{Func}(\mathcal{C},\mathcal{D})$ is the category
  whose objects are the functors from $\mathcal{C}$ to $\mathcal{D}$
  and whose \m s are the natural transformations.
\end{defn}

\begin{defn}
  Two categories, $\mathcal C$ and $\mathcal{D}$, are {\em
    isomorphic\/} ({\em equivalent\/}) when there are functors
  $\xymatrix@1{ {\mathcal{C}} \ar@<.5ex>[r]^F & {\mathcal{D}}
    \ar@<.5ex>[l]^G}$ such that the composite functors are (naturally
  isomorphic to) the respective identity functors.
%
\end{defn}

\begin{lemma}
  A functor \func{F}{\mathcal C}{\mathcal D} is an {\em equivalence\/}
  of categories if and only if
\begin{itemize}
\item any object of $\mathcal D$ is isomorphic to an object of the form
$F(a)$ for some object $a$ of $\mathcal C$
\item $F$ is bijective on \m\ sets: The maps $\mathcal C(a,b)
  \xrightarrow{f \to F(f)} \mathcal D(F(a),F(b))$ are bijections for
  all objects $a$ and $b$ of $\mathcal C$
\end{itemize}
\end{lemma}
\begin{proof}
  Suppose that \func{F}{\mathcal C}{\mathcal D} is an equivalence of
  categories. Then there is a functor $G$ in the other direction and
  natural iso\m s $\sigma \colon GF \Longrightarrow 1_{\mathcal{C}}$
  and $\tau \colon FG \Longrightarrow 1_{\mathcal{D}}$. Let $d$ be any
  object of $\mathcal{D}$. The iso\m\ $\tau_d \colon FG(d)
  \xrightarrow{\cong} d$ shows that $d$ is isomorphic to
  $Fa$ for $a=Gd$. Let $a,b$ be objects of $\mathcal{C}$. We note
  first that $\mathcal{C}(a,b) \to \mathcal{D}(Fa,Fb) \to
  \mathcal{C}(GFa,GFb)$ is injective for the commutative diagram
  \begin{equation*}
    \xymatrix{
      a \ar[d]_f & GFa \ar[d]_{GFf} \ar[r]^{\sigma_a}_{\cong} &
      a \ar[d]^f \\
      b & GFb \ar[r]^{\sigma_b}_{\cong} & b}
  \end{equation*}
  shows that $f=\sigma_b \circ GFf \circ \sigma_a^{-1}$ can be
  recovered from $GFf$. Thus $\mathcal{C}(a,b) \to \mathcal{D}(Fa,Fb)$
  is injective. Symmetrically, also the functor $G$ is injective on \m\
  sets. To show that $F$ is surjective on \m\ sets let $g$ be any
  $\mathcal{D}$-\m\ $Fa \to Fb$. Put $f=\sigma_b \circ Gg \circ
  \sigma_b^{-1}$. The commutative diagram
  \begin{equation*}
    \xymatrix{
      Fa \ar[d]_g &
      GFa \ar[d]_{Gg} \ar[r]^{\sigma_a}_{\cong} & a \ar[d]_f &
      GFa \ar[l]_{\sigma_a}^{\cong} \ar[d]^{GFf} & a \ar[d]^f \\
      Fb & 
      GFb \ar[r]^{\sigma_b}_{\cong} & b & GFb
      \ar[l]_{\sigma_b}^{\cong} & b }
  \end{equation*}
  shows that $GFf=Gg$ and so $Ff=g$ since $G$ is injective on \m\
  sets.
  
  Conversely, suppose that $F \colon \mathcal{C} \to \mathcal{D}$ is a
  functor satisfying the two conditions. We must construct a functor
  $G$ in the other direction and natural iso\m s $\tau \colon FG
  \Longrightarrow 1_{\mathcal{D}}$ and $\sigma \colon GF
  \Longrightarrow 1_{\mathcal{C}}$. By the first condition, for every
  object $d \in \mathcal{D}$, we can find an object $Gd \in
  \mathcal{C}$ and an iso\m\ $\tau_d \colon FGd \to d$. By the second
  condition, $\mathcal{C}(Gc,Gd) \cong \mathcal{D}(FGc,FGd)$ for any
  two objects $c$ and $d$ of $\mathcal{D}$. Here, $\mathcal{D}(c,d)
  \cong \mathcal{D}(FGc,FGd)$ because $FGc \cong c$ and $FGd \cong d$.
  Thus we have $\mathcal{D}(c,d) \cong \mathcal{D}(FGc,FGd) \cong
  \mathcal{C}(Gc,Gd)$. This means that for every $\mathcal{D}$-\m\
  \func{g}{c}{d} there is exactly one $\mathcal{C}$-\m\
  \func{Gg}{Gc}{Gd} such that
  \begin{equation*}
    \xymatrix{
      FGc \ar[d]_{FGg} \ar[r]^{\tau_c}_{\cong} & c \ar[d]^g \\
      FGd \ar[r]^{\tau_d}_{\cong} & d}
  \end{equation*}
  commutes. Now $G$ is a functor and $\tau$ a natural iso\m\ $FG
  \Longrightarrow 1_{\mathcal{D}}$. What about $GF$? Well, for any
  object $a$ of $\mathcal{C}$, $\mathcal{C}(GFa,a) \cong
  \mathcal{D}(FGFa,Fa) \ni \tau_{Fa}$ so there is a unique iso\m\
  $\sigma_a \colon GFa \to a$ such that $F\sigma_a=\tau_{Fa}$. This
  gives the natural iso\m\ $\sigma \colon GF \Longrightarrow
  1_{\mathcal{C}}$.
\end{proof}

It follows that when $\xymatrix@1{ {\mathcal{C}} \ar@<.5ex>[r]^F &
    {\mathcal{D}} \ar@<.5ex>[l]^G}$ is an equivalence of categories
  then there are bijections
\begin{equation*}
  \mathcal{C}(c,Gd) = \mathcal{D}(Fc,d) \qquad \qquad
  \mathcal{C}(Gd,c) = \mathcal{C}(d,Fc)
\end{equation*}
of \m\ sets.

\begin{lemma}\label{lemma:equivFunc}
  If $\mathcal{C},\mathcal{C}'$ and $\mathcal{D},\mathcal{D}'$ are
  equivalent, then the functor categories
  $\mathrm{Func}(\mathcal{C},\mathcal{D})$ and
  $\mathrm{Func}(\mathcal{C}',\mathcal{D}')$ are equivalent.
\end{lemma}

The {\em full subcategory\/} generated by some of the objects of
$\mathcal C$ is the category whose objects are these objects and whose
\m s are all \m s in $\mathcal C$. 

\begin{exmp}
  The category of finite sets is equivalent to the full subcategory
  generated by all sections $S_{<n}=\{x \in \Z_+ \mid x < n\}$, $n \in
  \Z_+$, of $\Z_+$. The category of finite dimensional real vector
  spaces is equivalent to the category $\mathbf{Mat}_{\mathbf R}$. If
  \func{f}{X}{Y} is a homeo\m\ (homotopy equivalence) then the induced
  \m\ $\pi(f) \colon \pi(X) \to \pi(Y)$ is an iso\m\ (equivalence) of
  categories. The fundamental groupoid of a space is
  equivalent to the full subcategory generated by a point in each path
  component.
\end{exmp}




\section{Categories of right $G$-sets}
\label{sec:actions}

Let $G$ be a topological group and $F$ and $Y$  topological spaces.

\begin{defn}\label{defn:Gaction}
  A {\em right\/} action of $G$ on $F$ is a continuous map $F \times G
  \to F \colon (x,g) \mapsto x \cdot g$, such that $x\cdot e=x$ and
  $x\cdot (gh)=(x \cdot g)\cdot h$ for all $g,h \in G$ and all $x \in
  F$.  A topological space equipped with a right $G$-action is called
  a right $G$-space.  A \co\ map \func{f}{F_1}{F_2} between two right
  $G$-spaces is a $G$-map if $f(xg)=f(x)g$ for all $g \in G$ and $x
  \in F_1$. 
\end{defn}

\begin{defn}\label{defn:leftGaction}
  A {\em left\/} action of $G$ on $Y$ is a continuous map $G \times Y
  \to Y \colon (g,y) \mapsto g \cdot y$, such that $e\cdot y=y$ and
  $(gh)\cdot y=g \cdot (h \cdot y) $ for all $g,h \in G$ and all $y \in
  Y$.  A topological space equipped with a left $G$-action is called
  a left $G$-space.  A \co\ map \func{f}{Y_1}{Y_2} between two left
  $G$-spaces is a $G$-map if $f(gy)=gf(y)$ for all $g \in G$ and $x
  \in Y_1$. 
\end{defn}

The orbit spaces (with the quotient topologies) are denoted $F / G =
\{ xG \mid x \in F \}$ for a right action $F \times G \to F$ and $G
\backslash Y =\{ Gy \mid y \in Y\}$ for a left action $G \times Y \to
Y$.

The {\em orbit\/} through the point $x \in F$ for the right action $F
\times G \to F$ is the sub-right $G$-space $xG=\{ xg \mid g \in G \}$
obtained by hitting $x$ with all elements of $G$; the {\em
  stabilizer\/} at $x$ is the subgroup ${}_xG = \{ g \in G \mid xg=x
\}$ of $G$. The universal property of quotient spaces gives a
commutative diagram
\begin{equation*}
    \xymatrix{
      G \ar[rr]^{g \to xg} \ar[dr]_{g \to {}_xGg} && xG \\
      & {}_xG \backslash G \ar[ur]_{ {}_xGg \to xg} }
\end{equation*}
of right $G$-spaces and $G$-maps
(\href{http://www.math.ku.dk/~moller/e03/3gt/notes/gtnotes.pdf}{General
  Topology, 2.81}).  
Note that $G$-map ${}_xG\backslash G \to xG
\colon {}_xGg \to xg$ is bijective. (In particular, the index of the
stabilizer subgroup at $x$ equals the cardinality of the orbit through
$x$.)  In many cases it is even a homeo\m\ so that the orbit $xG$
through $x$ and the coset space $_{x}G \backslash G$ of the isotropy
subgroup at $x$ are homeomorphic.

\begin{prop}[$G$-orbits as coset spaces]\label{prop:orbit}
  Suppose that $F$ is a right $G$-space and $x$ a point of
  $F$. Then
  \begin{equation*}
    \text{ ${}_xG \backslash G \to xG$ is a homeo\m} \iff
    \text{ $G \xrightarrow{g \to xg} xG $ is a quotient map}
  \end{equation*}
\end{prop}
\begin{proof}
  Use that the a bijective quotient map is a homeo\m , the composition
  of two quotient maps is quotient, and if the composition of two maps
  is quotient than the last map is quotient
  (\href{http://www.math.ku.dk/~moller/e03/3gt/notes/gtnotes.pdf}{General
  Topology, 2.77}). 
  By definition, $G \to {}_xG \backslash G$  is
  quotient.
\end{proof}

By a right (or left) $G$-{\em set\/} we just mean a right (or left)
$G$-space with the discrete topology. In the following we deal with
$G$-sets rather than $G$-spaces.

\begin{defn}\label{defn:GSet}
  $G\mathbf{Set}$ is the category of right $G$-sets and $G$-maps. The
  objects are right $G$-sets $F$ and the \m s $\varphi \colon F_1 \to
  F_2$ are $G$-maps (meaning that $\varphi(x g) = \varphi(x)g$ for all
  $x \in F_1$ and $g \in G$).
\end{defn}

\setcounter{subsection}{\value{thm}}
\subsection{Transitive right actions}
\label{sec:transitiveaction}
\addtocounter{thm}{1}

The right $G$-set $F$ is {\em transitive\/} if $F$ consists of a
single orbit.  If $F$ is transitive then $F=xG$ for some (hence any)
point $x \in F$ so that $F$ and $H \backslash G$ are isomorphic
$G$-sets where $H$ is the stabilizer subgroup at the point $x$
(Proposition~\ref{prop:orbit}). Thus any transitive right $G$-set is
isomorphic to the $G$-set $H \backslash G$ of right $H$-cosets for
some subgroup $H$ of $G$.


\begin{defn}\label{defn:OG}
  The orbit category of $G$ is the full subcategory $\mathcal O_G$ of
  $G\mathbf{Set}$ generated by all transitive right $G$-sets.
\end{defn}

 

The orbit category $\mathcal{O}_G$ of $G$ is equivalent to the full
subcategory of $G\mathbf{Set}$ generated by all $G$-sets of the form
$H \backslash G$ for subgroups $H$ of $G$. What are the \m s in the
orbit category $\mathcal{O}_G$?

\begin{defn}\label{defn:transp}
Let $H_1$ and $H_2$ be subgroups of $G$. The {\em transporter\/} is
the set 
\begin{equation*}
N_G(H_1,H_2)=\{n \in G \mid nH_1n^{-1} \subset H_2\}  
\end{equation*}
of group elements conjugating $H_1$ into $H_2$.  
\end{defn}

The transporter set $N_G(H_1,H_2)$ is a left $H_2$-set. Let
Let $H_2  \backslash N_G(H_1,H_2)$ be the set of $H_2$-orbits.


\begin{prop}\label{rightGmaps}
  There is a bijection
  \begin{equation*}
    \tau \colon H_2 \backslash N_G(H_1,H_2) \to 
    \mathcal O_G(H_1 \backslash G,H_2 \backslash G), \quad
    \tau({H_2n})(H_1g)=H_2ng  
  \end{equation*}
  This map takes $H_2n$ to left multiplication 
  $H_1 \backslash G \xrightarrow {H_1g \to H_2ng} H_2 \backslash G$
  by $H_2n$. In case $H_1=H=H_2$, the map
  \begin{equation*}
    \tau \colon H \backslash N_G(H) \to 
    \mathcal O_G(H \backslash G,H \backslash G), \quad
    \tau({Hn})(Hg)=Hng 
  \end{equation*}
  is a group iso\m .
\end{prop}
\begin{proof}
  The inverse to $\tau$ is the map that takes a
  $G$-map $H_1 \backslash
  G \xrightarrow{\varphi} H_2 \backslash G$ to its value
  $\varphi(H_1)=H_2n$ at $H_1 \in H_1 \backslash
  G$. Since $H_2n=\varphi(H_1)=\varphi(H_1H_1)=H_2nH_1$, the group
  element $n$ conjugates $H_1$ into $H_2$.
  In case $H_1=H=H_2$ and $n_1,n_2 \in N_G(H)$,
  we have
  \begin{equation*}
    \tau({Hn_1})\tau({Hn_2})(H) =
    \tau({Hn_1})(Hn_2)=Hn_1n_2=\tau({Hn_1n_2})(H) 
  \end{equation*}
  so $\tau$ is group homo\m\ in this case.
\end{proof}



In particular we see that 
\begin{itemize}
\item all \m s in $\mathcal{O}_G$ are epi\m s   
\item all endo\m s in $\mathcal O_G$ are auto\m s
\item every object $H \backslash G$ of $\mathcal{O}_G$ is equipped
  with left and right actions\set
  \begin{equation}
    \label{exmp:Gleftright}
    H \backslash N_G(H) \times  H \backslash
    G \times G =
    \mathcal{O}_G(H \backslash G, H \backslash G) \times H \backslash
    G \times G \to  H \backslash G \colon
    Hn \cdot  Hg \cdot m = Hngm
  \end{equation}\add
  where the left action are the $G$-auto\m s of $H \backslash G$ in
  $\mathcal{O}_G$.
\item the maximal $G$-orbit is $G=\{e\} \backslash G $ and
  $\mathcal{O}_G(\{e\} \backslash G, H \backslash G) = H \backslash
  G$, the minimal $G$-orbit is $\ast = G \backslash G$ and
  $\mathcal{O}_G(H \backslash G, G \backslash G) = \ast$ ($G
  \backslash G = \ast$ is the final object of $\mathcal{O}_G$)
\end{itemize}


\begin{rmk}[Iso\m\ classes of objects of $\mathcal{O}_G$]
  The set of objects of $\mathcal{O}_G$ corresponds to the set of
  subgroups of $G$. The set of iso\m\ classes of objects of
  $\mathcal{O}_G$ corresponds to the set of conjugacy classes of
  subgroups of $G$: Two objects $H_1 \backslash G$ and $H_2 \backslash
  G$ of the orbit category $\mathcal O_G$ are isomorphic if and only
  if $H_1$ and $H_2$ are conjugate: If there exist an inner auto\m\ 
  that takes $H_1$ into $H_2$ and an inner auto\m\ that takes $H_2$
  into $H_1$ such that the composite maps are the respective identity
  maps of $H_1 \backslash G$ and $H_2 \backslash G$, then these inner
  auto\m s must in fact give bijections between $H_1$ and $H_2$ as the
  factorizations $H_1 \xrightarrow{\mathrm{Inn}(n_1)} H_2
  \xrightarrow{\mathrm{Inn}(n_2)} H_1 \xrightarrow{\mathrm{Inn}(n_1)}
  H_2$ of the respective identity maps imply that the inner auto\m\ 
  $\mathrm{Inn}(n_1)$ is a bijection.
\end{rmk}



\section{The classification theorem}
\label{sec:1conncov}

In this section we shall see that covering maps are determined by
their monodromy functor.

\begin{defn}
  $\mathrm{Cov}(X)$ is the category of covering spaces over the space
  $X$. The objects are covering maps $Y \to X$ and the \m s
  $\mathrm{Cov}(X)(p_1 \colon Y_1 \to X, p_2 \colon Y_2 \to X)$ are
  \co\ maps \func f{Y_1}{Y_2} over $X$ (meaning that $f$ preserves
  fibres or $p_1=p_2f$).
\end{defn}

How can we describe the category $\mathrm{Cov}(X)$? We are going to
assume from now on that $X$ is {\em path connected and locally path
connected\/}. 

Let $\mathrm{Func}(\pi(X),\mathbf{Set})$ be the category of functors
from the fundamental groupoid $\pi(X)$ to the category $\mathbf{Set}$
of sets. There is a functor
\begin{equation*}
  \mathrm{Cov}(X) \to \mathrm{Func}(\pi(X),\mathbf{Set})
\end{equation*}
which takes a covering map \func pYX to its monodromy functor $\func
{F(p)}{\pi(X)}{\mathbf{Set}}$ \eqref{eq:defnp-1} and a covering map
\m\ to the induced natural transformation of functors.  Conversely,
does any such functor come from a covering space of $X$?

Suppose that $\func F{\pi(X)}{\mathbf{Set}}$ is any functor. Let $Y(F)
= \bigcup_{x \in X}F(x)$ be the union of the fibres and let \func
{p(F)}{Y(F)}X be the obvious map taking $F(x)$ to $x$ for any point $x
\in X$.

\begin{defn}
  A space $X$ is semi-locally simply connected at the point $x \in X$ 
   if any \n\ of $x$ contains a \n\ $U$ of $x$ such that any loop at
   $x$ in $U$ is contractible in $X$. The space $X$ is semi-locally
   simply connected if it is  semi-locally simply connected at all its
   points.
\end{defn}

All locally simply connected spaces are semi-locally simply connected.

\begin{lemma}
  Suppose that $X$ is locally path connected and semi-locally simply
  connected. Then there is a topology on $Y(F)$ such that
  \func{p(F)}{Y(F)}X is a covering map. The monodromy functor of
  \func{p(F)}{Y(F)}X is $F$.
\end{lemma}
\begin{proof}
  Suppose that $x$ is a point in $X$ and $U \subset X$ an open path
  connected \n\ of $x$ such that any loop in $U$ based at $x$ is
  nullhomotopic in $X$. Observe that this implies that there is a
  unique path homotopy class $u_z$ from $x$ to any other point $z$ in
  $U$ so that
   \begin{equation*}
  U \times F(x) \to p^{-1}(U) \colon (y,z) \to F(u_z)(y)
\end{equation*}
is a bijection. 

For each $y \in F(x)$, let $(U,y) \subset Y$ be the image of $U \times
\{y\}$ under the above bijection.  By assumption, the topological
space $X$ has a basis of sets $U$ as above. The sets $(U,y)$ then form
a basis for a topology on $Y$.

The covering map $Y(F) \to X$ determines a fibre functor
\eqref{eq:defnp-1} from the fundamental groupoid of $X$ to the
category of sets. By construction, this fibre functor is $F$.
\end{proof}

\begin{defn}  \label{defn:X1} 
  A covering map \func pYX is  universal if $Y$ is simply connected.  
\end{defn}

According to the Lifting Theorem~\ref{thm:lift}, any two universal
covering spaces over $X$ are isomorphic in the category
$\mathrm{Cov}(X)$ of covering spaces over $X$. We may therefore speak
about {\em the\/} universal covering space of $X$. Is there always a
universal covering space of $X$?

By Corollary~\ref{cor:pi1Yinjects} the fundamental groupoid of $Y(F)$
has the set $Y(F)$ as object set and the \m s are
\set\begin{equation}\label{eq:piYF} \pi(Y(F))(y_1,y_2) = \{ u \in
  \pi(X)(x_1,x_2) \mid F(u)y_1 = y_2 \}
\end{equation}\add
for all points $x_1,x_2 \in X$ and $y_1 \in F(x_1)$, $y_2 \in F(x_2)$.
In particular, let $x_0$ be a base point in $X$. There is a right
action $F(x_0) \times \pi_1(X,x_0) \to F(x_0)$ and
\begin{align*}
  \text{$Y(F)$ is path connected} &\iff 
  \text{The right action of $\pi_1(X,x_0)$ on $F(x_0)$ is 
    transitive} \\
   \text{$Y(F)$ is simply connected} &\iff 
  \text{The right action of $\pi_1(X,x_0)$ on $F(x_0)$ is simply transitive}
\end{align*}
We can always find a functor that satisfies the last condition in that
\begin{equation*}
 F= \pi(X)(x_0,-) \colon \pi(X) \to \mathbf{Set}
\end{equation*}
is a functor and the action of $\pi_1(X,x_0)$ on $F(x_0)=\pi_1(X,x_0)$
is simply transitive.

\begin{cor}
  $X$  admits a simply
  connected covering space if and only if $X$ is semi-locally simply
  connected. 
\end{cor}
\begin{proof}
  The covering space $Y(F)$ of the functor $ F= \pi(X)(x_0,-)$ is
  simply connected.

Conversely, suppose that \func pYX is a covering map and $U \subset X$
and evenly covered open subspace then $U \to X$ factors through $Y \to
X$. If $\pi_1(Y)$ is trivial then $\pi_1(U) \to \pi_1(X)$ is the
trivial homo\m .
\end{proof}

\begin{exmp}
  The Hawaiian Earring $\bigcup_{n \in \Z_+}C_{1/n}$ and the infinite
  product $\prod S^1$ of circles are connected and locally
  path connected but not semi-locally simply connected.  Thus they
  have no simply connected covering spaces. The infinite join $\bigvee
  S^1$ does have a simply connected covering space since it is a
  CW-complex. Indeed any CW-complex or manifold is locally
  contractible \cite[Appendix]{hatcher}, in particular locally simply
  connected.
\end{exmp}

\begin{thm}[Classification of Covering Maps]\label{thm:classcover2}
  Suppose that $X$ is semi-locally simply connected.  The monodoromy
  functor and the functor $F \to Y(F)$
  \begin{equation*}
    \xymatrix@1{
      {\mathrm{Cov}(X)} \ar@<.5ex>[r] &
      {\mathrm{Func}(\pi(X),\mathbf{Set})} \ar@<.5ex>[l] }
  \end{equation*}
  are category iso\m s.
\end{thm}
\begin{proof}
Let \func {p_1}{Y_1}X and  \func {p_2}{E_Y}X be covering maps over $X$
with associated functors $F_1$ and $F_2$. A covering map
\begin{equation*}
  \xymatrix{ Y_1 \ar[dr]_{p_1}  \ar[rr]^f && 
  Y_2 \ar[dl]^{p_2} \\ & X}
\end{equation*}
induces a natural transformation $\tau_f \colon F_1 \Longrightarrow F_2$ of
functors given by $\tau_f(x)=f \vert p_1^{-1}x \colon p^{-1}x_1 \to
p^{-1}x_2$.  Conversely, any natural transformation $\tau \colon F_1
\Longrightarrow F_2$ induces a covering map $Y(f) \colon Y(F_1) \to
Y(F_2)$ of the associated covering spaces.
\end{proof}

For example, let $\func F{\pi(X)}{\mathbf{Set}}$ be any functor, let
$x_0 \in X$ and $y_0 \in F(x_0)$.  There is a natural transformation
$\pi(X)(x_0,-) \Longrightarrow F$ whose $x$-component is
$\pi(X)(x_0,x) \to F(x) \colon u \to F(u)y_y$ for any point $x$ of
$X$. This confirms that the universal covering space lies above them
all.

\begin{cor}
  The functor
  \begin{equation*}
    \mathrm{Cov}(X) \to \pi_1(X,x_0)\mathbf{Set} \colon 
 (\func pYX) \to p^{-1}(x_0)
  \end{equation*}
  is an equivalence of categories.
\end{cor}
\begin{proof}
  The inclusion $\pi_1(X,x_0) \to \pi(X)$ 
  of the the full subcategory of $\pi(X)$ generated by
  $x_0$ into $\pi(X)$ is an equivalence of categories. The induced
  functor ${\mathrm{Func}(\pi(X),\mathbf{Set})} \to
  {\mathrm{Func}(\pi_1(X,x_0),\mathbf{Set})}$ is then also an
  equivalence. But ${\mathrm{Func}(\pi_1(X,x_0),\mathbf{Set})}$ is
  simply the category of right $\pi_1(X,x_0)$-sets.
\end{proof}

In particular, the full subcategory $\mathrm{Cov}_0(X)$ of connected
covering spaces over $X$ is equivalent to the category of transitive
right ${\pi_1(X,x_0)}$-sets which again is equivalent to the orbit
category $\mathcal{O}_{\pi_1(X,x_0)}$ (\ref{defn:OG}). The set of
covering space \m s from the connected covering space \func{p_1}{Y_1}X
to the connected covering space \func{p_2}{Y_2}X is
\begin{align*}
  \mathrm{Cov}(X)(\func{p_1}{Y_1}X,\func{p_2}{Y_2}X) &=
  \mathrm{Func}(\pi(X),\mathbf{Set})(F(p_1),F(p_2)) \\
  &= \pi_1(X)\mathbf{Set}(p_1^{-1}(x_0),p_2^{-1}(x_0)) \\
  &= \mathcal{O}_{\pi_1(X)}(\pi_1(Y_1) \backslash \pi_1(X), 
  \pi_1(Y_2) \backslash \pi_1(X)) \\
  &= \pi_1(Y_2) \backslash N_{\pi_1(X)}(\pi_1(Y_1),\pi_2(Y_2))
\end{align*}
and, in particular,
\begin{equation*}
   \mathrm{Cov}(X)(\func pYX,\func pYX) = 
   \pi_1(Y) \backslash N_{\pi_1(X)}(\pi_1(Y))
\end{equation*}
for any connected covering space \func pYX over $X$. If we map out of
the universal covering space $X \langle 1 \rangle \to X$ this gives
\begin{equation*}
   \mathrm{Cov}(X)(X \langle 1 \rangle \to X, Y \to X) = \pi_1(Y)
   \backslash \pi_1(X) \qquad\qquad
    \mathrm{Cov}(X)(X \langle 1 \rangle \to X, X \langle 1 \rangle \to
    X) = \pi_1(X)
\end{equation*}
which means that the universal covering space admits a left covering
space $\pi_1(X)$-action with orbit space $\pi_1(X) \backslash X
\langle 1 \rangle \to X = X$.

\begin{cor}
  Let $G=\pi_1(X)$ for short. The functor
  \begin{equation*}
    \mathcal{O}_{G} \to \mathrm{Cov}_0(X) \colon H \to 
    (H  \backslash X \langle 1 \rangle \to  G \backslash X \langle 1 \rangle)
   \end{equation*}
   is an equivalence of categories.
\end{cor}

  Is this
  \begin{equation*}
    \xymatrix{ 
    & C_2 \backslash \Sigma_3 \ar[dr] \ar@(ur,ul)[]|{\mathrm{id}} \\
     *\txt{\makebox[5mm][l]{$\{e\} \backslash \Sigma_3$}}
       \ar@(dl,ul)[]|{\Sigma_3}
           \ar@<-.5ex>[ur] \ar[ur]
      \ar@<.5ex>[ur] \ar@<-.5ex>[dr] \ar@<.5ex>[dr]  && 
      *\txt{\makebox[5mm][r]{$\Sigma_3 \backslash \Sigma_3$}} 
     \ar@(dr,ur)[]|{\mathrm{id}} \\
      &  {C_3 \backslash \Sigma_3} \ar[ur] \ar@(dr,dl)[]|{C_2}}
  \end{equation*}
  a picture of the orbit category of symmetric group $\Sigma_3$ or is
  it a picture of the path connected covering spaces over a path
  connected, locally path connected, and semi-locally simply connected
  space with fundamental group $\Sigma_3$? Both! The space could be
  $X_{\Sigma_3}$ from Corollary~\ref{cor:XG}; see
  Example~\ref{exmp:RP2astRP2} for more information.

Here are some examples to illustrate the Classification of Covering
Spaces. 
\begin{description}
\item[Covering spaces of the circle] The category
  $\mathrm{Cov}_0(S^1)=\mathcal{O}_{C_\infty}$ of path connected covering
  spaces of the circle $S^1 = \Z \backslash \R$ consists of the
  covering spaces $n\Z \backslash \R \to \Z \backslash \R$ where
  $n=0,1,2,\ldots $. There is a covering map $n\Z \backslash \R \to
  m\Z \backslash \R$ if and only if $m \vert n$ and in that case there
  are $m$ such covering maps, namely the maps
  \begin{equation*}
    \xymatrix{
             S^1 \ar[dr]_{z^n} \ar[rr]^{\zeta z^{m/n}} && 
             S^1 \ar[dl]^{z^m} \\
             & S^1 }
  \end{equation*}
  where $\zeta$ is any $m$th root of unity.
\item[Covering spaces of projective spaces] The category
  $\mathrm{Cov}_0(\R P^n) = \mathcal{O}_{C_2}$ of connected covering
  spaces of real projective $n$-space $\R P^n$, $n \geq 2$, has $2$
  objects, namely the trivial covering map $\R P^n \to \R P^n$ and the
  universal covering map $S^n \to \R P^n$.
\item[Covering spaces of lense spaces]
  The universal covering space of
  the lense space $L^{2n+1}(m)=C_m\backslash S^{2n+1}$, $n \geq 1$, is
  $S^{2n+1}$. The other covering spaces are the lense spaces
  $L^{2n+1}(r)=C_r\backslash S^{2n+1}$ for eah divisor $r$ of $m$. The
  category of connected covering spaces of $L^{2n+1}(m)$ is equivalent
  to the orbit category $\mathcal{O}_{C_m}$.
\item[Covering spaces of surfaces] The category
  $\mathrm{Cov}_0(M_g)=\mathcal{O}_{\pi_1(M_g)}$ is harder to describe
  explicitly.  Any finite sheeted covering space of a compact surface
  is again a compact surface.  The paper \cite{mednykh88} contains
  information about covering spaces of closed surfaces.
\end{description}


\begin{exmp}[Covering spaces of the M{\"o}bius band]
  The cylinder $S¹\times [-1,1] = \Z \backslash (\R \times
  [-1,1])$ where the action is given by $n \cdot (x,t) \to (x+n,t)$

\noindent \hfill
\begin{picture}(10,60)(410,-5)
\thicklines
\put(0,50){\vector(1,0){400}}
\put(0,0){\vector(1,0){400}}
\multiput(10,0)(50,0){8}{\vector(0,1){29}}
\multiput(10,0)(50,0){8}{\line(0,1){50}}
\multiput(20,40)(50,0){8}{$\bullet$}
\end{picture}

\noindent
The M{\"o}bius band $\mathrm{MB} =  \Z \backslash (\R \times
  [-1,1])$ where the action is given by $n \cdot (x,t) \to (x+n,(-1)^nt)$

\noindent \hfill
\begin{picture}(10,60)(410,-5)           
\thicklines
\put(0,50){\vector(1,0){400}}
\put(0,0){\vector(1,0){400}}
\multiput(10,50)(100,0){4}{\vector(0,-1){29}}
\multiput(10,50)(100,0){4}{\line(0,-1){50}}
\multiput(60,0)(100,0){4}{\vector(0,1){29}}
\multiput(60,0)(100,0){4}{\line(0,1){50}}
\multiput(20,40)(100,0){4}{$\bullet$}
\multiput(70,10)(100,0){4}{$\bullet$}
\end{picture}

\noindent
Every even-sheeted covering space of the M{\"o}buis band is a
cylinder, every odd-sheeted covering space is a  M{\"o}buis band. 

\noindent \hfill
\begin{picture}(0,0)(-20,-10)
\thicklines
\put(0,0){\circle{50}}
\put(0,0){\circle{10}}
\put(0,5){\line(0,1){15}}
\put(0,-5){\line(0,-1){15}}
\put(5,0){\line(1,0){15}}
\put(-5,0){\line(-1,0){15}}
\end{picture}
\end{exmp}

\begin{exmp}[Covering spaces of $S¹ \cup_m (S¹ \times I)$]
  Let $X_m=S¹ \cup_m (S¹ \times I)$ be the mapping cylinder of the
  degree $m$ map of the circle. We can construct $X_m$ in the
  following way: Take a (codomain) circle of circumference $1/m$ and a
  square $[0,1] \times [0,1]$. Wrap the bottom edge $[0,1] \times
  \{0\}$ of the square $m$ times around the circle in a screw motion
  so that each time the square travels once around the circle it is
  also being rotated an angle of $2\pi/m$.  Finally, glue the two
  ends, $\{0\} \times [0,1]$ and $\{1\} \times [0,1]$, of the square
  together.  There is a picture of $X_m$ in \cite[Example
  1.29]{hatcher}.  The codomain circle is the core circle and the
  domain circle is the boundary circle. The fundamental group
  $\pi_1(X_m)$ is $\Z$ since $X_m$ deformation retracts onto the
  codomain (core) circle so that the inclusion $S¹ \xrightarrow{i_1}
  X_m \supset S¹$ is a homotopy equivalence.  The inclusion $S¹
  \xrightarrow{i_0} X_m \supset S¹$ of the domain (boundary) circle
  induces multiplication by $m$ on the fundamental groups; this is
  simply because of the general mapping cylinder diagram which becomes
  \begin{equation*}
    \xymatrix{ S^1 \ar[dr]_m \ar[r]^-{i_0} & (S¹ \times I) \cup_m S¹=X_m
      \ar@<1.2ex>[d]^{\mathrm{DR}} \\
      & S¹ \ar@<1.2ex>[u]^{i_1}_{\simeq}}
  \end{equation*}
  in this special case. It may help to envision the boundary circle in
  $X_m$ sliding towards the core circle.
 
  The universal covering space of $X_m$ is $X_m \langle 1 \rangle=
  C\Z/m \times \R$ where $C\Z/m=(\Z/m \times I) /(\Z/m \times \{1\})$
  is the cone on the set $\Z/m$ with $m$ points. ($C\Z/m$ is a
  starfish with $m$ arms). 
  \begin{figure}[h]\label{fig:ex6_7}
  \centering
   \includegraphics[width=5cm,height=2cm]{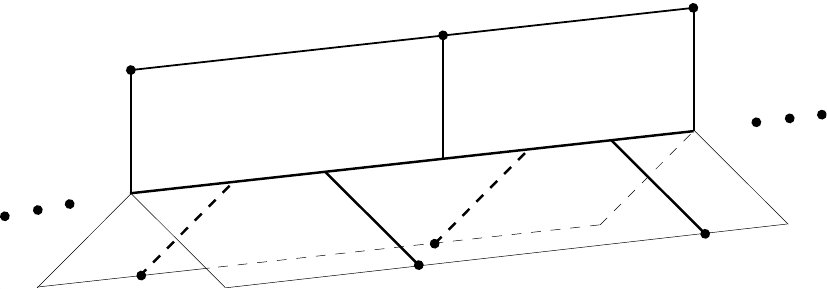}
  \caption{The universal covering space of $X_m$}
  \end{figure} 
  We may realize $C\Z/m \times \R$ in $\R³$
  with $C\Z/m$ placed horizontally in the $XY$-plane and $\R$ as the
  vertical $Z$-axis.  
 The covering space action of the unit $1 \in \Z$
  on $C\Z/m \times \R$ is then the screw motion $([[a]_m,t],x) \to
  ([[a+1]_m,t],x+1/m)$ with matrix
  \begin{equation*}
    \begin{pmatrix}
      \cos(2\pi/m) & -\sin(2\pi/m) & 0 \\
      \sin(2\pi/m) & \phantom{-}\cos(2\pi/m) & 0 \\
      0 & 0 & 1/m
    \end{pmatrix}
  \end{equation*}
  that rotates $C\Z/m$ counterclockwise $1/m$th of a full rotation and
  moves up along the $Z$-axis $1/m$th of a unit. 
  (In Figure~\ref{fig:ex6_7} the $\R$-axis isn't
  exactly vertical since that would take up too much space. The
  covering space action takes the indicated lines, situated at
  distance $1/m$, to each other.)
  What is the lift of
  the domain and the codomain circles of $X_m$ to the universal
  covering space $X_m \langle 1 \rangle$? (One of them will lift to a
  loop.)

  Since $m \in \Z$ acts trivially on $C\Z/m$ there is an $m$-sheeted
  covering map 
  \begin{equation*}
  C\Z/m \times S¹= C\Z/m \times m\Z \backslash \R =
                  m\Z \backslash (C\Z/m \times \R) \to 
                   \Z \backslash (C\Z/m \times \R) = X_m 
  \end{equation*}
  with $m\Z \backslash \Z$ as deck transformation group. What is the
  lift of the domain and the codomain circles to this $m$-fold
  covering space?
\end{exmp}

Let $X=X_1 \cup X_2$ be a CW-complex that is the union of two
connected subcomplexes $X_1$ and $X_2$ with connected intersection
$X_1 \cap X_2$. According to van Kampen, the fundamental group
$G=\pi_1(X)=G_1 \amalg_A G_2$ is the free product of $G_1=\pi_1(X_1)$ and
$G_2=\pi_1(X_2)$ with $A=\pi_1(X_1 \cap X_2)$ amalgamated. We will
assume that the homo\m s $G_1 \leftarrow A \rightarrow G_2$ are
injective. Then also the homo\m s $G_1 \rightarrow G \leftarrow G_2$ of
the push-out diagram
\begin{equation*}
  \xymatrix{
             A \ar[d] \ar[r] & G_2 \ar[d] \\
             G_1 \ar[r] & G }
\end{equation*}
are injective according to the Normal Form Theorem for Free Products
with Amalgamation \cite[Thm 2.6]{LyndonSchupp77}.

Let $X \langle 1 \rangle$ be the universal covering space of $X=G
\backslash X \langle 1 \rangle$ and let \func{p}{X \langle 1
  \rangle}{X} be the covering projection map. The spaces $p^{-1}(X_1)$
and $p^{-1}(X_2)$ are left $G$-spaces with intersection $p^{-1}(X_1)
\cap p^{-1}(X_1)=p^{-1}(X_1 \cap X_2)$. Let $y_0 \in p^{-1}(X_1 \cap
X_2)$ be a base point. The commutative diagram \cite[II.7.5]{brown82}
\begin{equation*}
  \xymatrix{
   {\pi_1(p^{-1}X_1,y_0)} \ar@{^(->}[d] \ar[r] &
   {\pi_1(p^{-1}X,y_0)=\{1\}}  \ar@{^(->}[d] \\
   {\pi_1(X_1,p(y_0))} \ar@{^(->}[r] &
   {\pi_1(X,p(y_0))} }
\end{equation*}
tells us that the component of $p^{-1}(X_1)$ containing $y_0$ is
simply connected so it is the universal covering space $X_1 \langle 1
\rangle$ of $X_1=G_1 \backslash X_1 \langle 1 \rangle$. We see from
this that there is a homeo\m\ of left $G$-spaces
\begin{equation*}
 G \times_{G_1} X_1 \langle 1 \rangle  \stackrel{\simeq}{\to} p^{-1}(X_1)
\end{equation*}
induced by the map $G \times X_1 \langle 1 \rangle \to p^{-1}(X_1)$
sending $(g,y)$ to $gy$. Similar arguments apply to $p^{-1}(X_2)$ and
$p^{-1}(X_1 \cap X_2)$, of course, and hence
\begin{equation*}
  X \langle 1 \rangle =  G \times_{G_1} X_1 \langle 1 \rangle  
  \cup_{G \times_A (X_1 \cap X_2) \langle 1 \rangle }
   G \times_{G_2}X_2 \langle 1 \rangle 
\end{equation*}
is the union of the two $G$-spaces $G \times_{G_i} X_i \langle 1
\rangle$, $i=1,2$. This means that the universal covering space of $X$
is the union of the $G$-translates of the universal covering spaces of
$X_1$ and $X_2$ joined along $G$-translates of the universal covering
space of $X_1 \cap X_2$.  The next example demonstrates this principle.

\begin{exmp}\cite[1.24, 1.29, 1.35, 1.44, 3.45] {hatcher}
  Let $X_{mn} = X_m \cup_{S¹}X_n$ be the double mapping cylinder for
  the degree $m$ map and the degree $n$ map on the circle. $X_{mn}$ is
  the union of the two mapping cylinders with their domain (boundary)
  circles identified, $X_m \cap X_n = S^1$.  By van Kampen, the
  fundamental group has a presentation
  \begin{equation*}
    \pi_1(X_{mn})=\pi_1(X_m) \amalg_{\pi_1(S¹)} \pi_1(X_n) = \gen{a,b
      \mid a^m=b^n}=G_{mn}
  \end{equation*}
  with two generators and one relation. We shall now try to build its
  universal covering space. 
  
  We may equip $X_{mn}$ with the structure of a $2$-dimensional
  CW-complex.  The $1$-skeleton of $X_{mn}$ consists of two circles,
  $a$ and $b$, joined by an interval, $c$, and 
 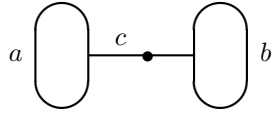
\begin{figure}[h]
  \centering
  \begin{picture}(75,40)
  \thicklines
  \multiput(0,20)(60,0){2}{\oval(20,40)}
  \put(10,20){\line(1,0){40}}
  \put(-20,18){$a$}
  \put(75,18){$b$}
  \put(20,24){$c$}
  \put(30,17){$\bullet$}
  \end{picture}
  \caption{$1$-skeleton $X^1_{mn}$ of $X_{mn}$}
  \label{fig:1xmn}
\end{figure}
$X_{mn}=X^1_{mn} \cup_{a^mc\overline{b}^n\overline{c}}D²$ is obtained
by attaching a $2$-cell along the loop
$a^mc\overline{b}^n\overline{c}$. (If we use the corollary to van Kampen
\cite[1.26]{hatcher} instead of the van Kampen theorem itself we get
that $\pi_1(X_{mn})=\gen{a,cb\overline{c} \mid
  a^m(c{b}\overline{c})^{-n}}$.)

The universal covering space $X_{mn} \langle 1 \rangle$ is also a
$2$-dimensional CW-complex.   The inverse image in $X_m \langle 1
\rangle$ of the left half of the $1$-skeleton is the vertical line
$\R$ with spiraling \lq rungs\rq\ attached $1/m$th of a unit apart.
Rungs with vertical distance $1$ point in the same direction so they
can be joined up with the inverse image in $X_n\langle 1 \rangle$ of
the right half of the $1$-skeleton.
  \begin{figure}[h]
  \centering
  \begin{picture}(75,75)
  \multiput(0,0)(60,0){2}{\line(0,1){75}}
  \multiput(0,10)(0,60){2}{\line(1,0){25}}
  \multiput(60,10)(0,60){2}{\line(-1,0){25}}
  \put(-15,35){$a^m$}
  \put(65,35){$b^n$}
  \multiput(22,8)(0,59){2}{$\bullet$}
  \multiput(34,8)(0,59){2}{$\bullet$}
  \end{picture}
  \caption{Part of $1$-skeleton of $X_{mn}\langle 1 \rangle$}
  \label{fig:XMN}
\end{figure}
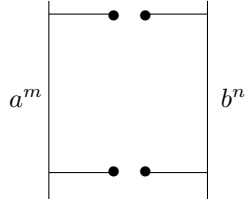
Now fill in $2$-cells in each of the rectangles with sides $a^m$, $c$,
$b^n$ and $c$.  Continue this process. There will be similar
rectangles shifted up $1/m$th unit along the left axis and rotated
$2\pi/m$ or up $1/n$th unit along the right axis and rotated $2\pi/n$.
The $2$-dimensional CW-complex $X_{mn} \langle 1 \rangle$ built in
this way is the universal covering space; it is the product $T_{mn}
\times \R$ of a tree $T_{mn}$ and the real line, hence contractible
\cite[Chp 3, Sec 7, Lemma 1]{spanier}.  The element $a \in G_{mn}$
acts by skew motion around one of the vertical lines in $X_m \langle 1
\rangle$ and $b \in G_{mn}$ acts by skew motion around one of the
vertical lines in $X_n \langle 1 \rangle$.  Note that $a^m=b^n$ acts
by translating one unit up. What is the lift of $X_m \cap X_n$ (the
circle parallel to circle $a$ but passing through the point $\bullet$
of the $1$-skeleton) to the universal covering space?

What is the universal abelian covering space $G_{mn}' \backslash
X_{mn} \langle 1 \rangle$ of $X_{mn}$? Its deck transformation group
is
\begin{equation*}
  G_{mn}' \backslash G_{mn} = (G_{mn})_{\mathrm{ab}}= 
  \gen{a,b|a^m=b^n,ab=ba} = \Z \times \Z/d
\end{equation*}
where $d=(m,n)$ is the greatest common divisor.  What is the $mn$ fold
covering space with fundamental group equal to the normal closure $N$
of $\gen{a^m,aba^{-1}b^{-1}}$ and deck transformation group $N
\backslash G = \gen{a,b \mid a^m=b^n,a^m,ab=ba} = \gen{a,b \mid
  a^m,b^n,ab=ba} = \Z/m \times \Z/n$? What is the lift of $X_m \cap
X_n$ to this covering space?
\end{exmp}

\setcounter{subsection}{\value{thm}}
\subsection{Cayley tables, Cayley graphs, and Cayley complexes} 
\cite[III.4]{LyndonSchupp77} \cite{brown82} 
\label{sec:cayley}\add
For any group presentation $G = \gen{g_{\alpha} \mid r_{\beta}}$ there
exists  (Corollary~\ref{cor:XG}) a $2$-dimensional CW-complex
\begin{equation*}
  X_{G \backslash G} = 
  D^0 \cup \coprod_{\{g_{\alpha}\}} D^1 \cup \coprod_{\{r_{\beta}\}}
  D^2 =
 (G \backslash G \times D^0) \cup 
 \coprod_{\{g_{\alpha}\}}(G \backslash G \times D^1) \cup 
 \coprod_{\{r_{\beta}\}}(G \backslash G \times D^2)  
\end{equation*}
with fundamental group $\pi_1(X_{G \backslash G})= \gen{g_{\alpha}
  \mid r_{\beta}} = G$. This is the most simple space with fundamental
group $G$ so it is natural to apply Theorem~\ref{thm:classcover2} to
$X_{G \backslash G}$.  So what are the connected covering spaces of
$X_G$? There is an equivalence of categories
\begin{equation*}
 {X}_? \colon \mathcal{O}_G \to \mathrm{Cov}_0(X_{G \backslash G}) \colon H
 \backslash G \to (X_{H \backslash G} \to X_{G \backslash G})
\end{equation*}
and the {\em Cayley complex\/} of $H \backslash G$ is the
$2$-dimensional CW-complex $X_{H\backslash G}$ while the {\em Cayley
  graph\/} is its $1$-skeleton.
%
We now define these CW-complexes more explicitly for any object of
$\mathcal{O}_G$ (or for any right $G$-space for that matter) relative
to the given presentation of $G$.
  
The $0$-skeleton of $X_{H \backslash G}$ is the right $G$-set $X^0_{H
  \backslash G} = H \backslash G$; this is the fibre
  of the covering map $X_{H \backslash G} \to X_{H \backslash G}$ as a
  right $G$-space. The $1$-skeleton of $X_{H
  \backslash G}$ is the {\em Cayley graph\/} for $H \backslash G$, the
$1$-dimensional $H \backslash N_G(H)$-CW-complex
  \begin{equation*}
  {X}^1_{H \backslash G} =
  (H \backslash G \times D^0)  \cup \coprod_{ \{Hg \to Hgg_{\alpha}\}} (H
  \backslash G \times D^1)
  \end{equation*}
  obtained from the $0$-skeleton $H \backslash G$ by attaching to each
  right coset $Hg \in H \backslash G$ an arrow from $Hg$ to
  $Hgg_{\alpha}$ for each generator $g_{\alpha}$; note that we have no
  other choice since the loop $g_\alpha$ in the base space lifts to a
  path in the total space that goes from $Hg$ in the fibre $H
  \backslash G$ to $Hgg_\alpha$ in the fibre. (The Cayley graph is
  simply a graphical presentation of the Cayles table for
  \href{http://en.wikipedia.org/wiki/Cayley_table} {Cayley table} for
  group multiplication $H \backslash G \times G \to H \backslash G$.)
  In this way, the Cayley table for $H
  \backslash G$ is a $\order{G \colon H}$-fold covering space of the
  $1$-skeleton $\bigvee_{\{g_{\alpha}\}}S^1$ of $X_G$.  The Cayley
  graph is connected since each group element $g$ is a product of the
  generators which means that there is a sequence of arrows connecting
  the $0$-cells $He$ and $Hg$.

  Next attach $2$-cells at each $Hg \in H \backslash G$ along the loop
  $r_{\beta}$ for each relation $r_{\beta}$. Since the relation
  $r_{\beta}$ is a factorization of the neutral element $e$ in terms
  of the $g_{\alpha}$, it defines loops $Hg \to Hgr_{\beta}=Hg$ based
  at each $0$-cell $Hg$ in the Cayley graph ${X}^1_{H \backslash G}$.
  The resulting left $H \backslash N_G(H)$-CW-complex
  \begin{equation*}
  {X}_{H \backslash G}=
   (H \backslash G \times D^0)  \cup 
   \coprod_{ \{Hg \to Hgg_{\alpha}\}} (H \backslash G \times D^1)  \cup 
   \coprod_{ \{Hg \xrightarrow[r_\beta]{} Hg\}} (H \backslash G \times D^2)
  \end{equation*}
  is the {\em Cayley complex\/} of $H \backslash G$. The Cayley
  complex is still connected for attaching $2$-cells does not alter
  the set of path components (Corollary~\ref{cor:pi1CW}). Clearly,
  every $G$-map $H_1 \backslash G \to H_2 \backslash G$ extends to a
  covering map $ {X}_{H_1 \backslash G} to  {X}_{H_2 \backslash G}$.

  In particular, taking $H=\{e\}$ to be the trivial group, the Cayley
  complex for the right $G$-set $\{e\} \backslash G = G$, 
  \begin{equation*}
    {X}_{\{e\} \backslash G}=
   (G \times D^0) \cup  \coprod_{\{g_{\alpha}\}} (G \times D^1) \cup
   \coprod_{ \{r_{\beta}\}} (G \times D^2) 
  \end{equation*}
  is a $2$-dimensional left $G$-CW-complex, the universal covering
  space of $X_{G \backslash G}$.  The $0$-skeleton is $G$, at
  each $g \in G$ there is an arrow from $g$ to $gg_{\alpha}$ for each
  generator $g_{\alpha}$ and a $2$-cell attached by the loop $g
  \xrightarrow{r_{\beta}} gr_{\beta}=g$. In other words, there is one
  $0$-$G$-cell $G \times D^0$, one $G$-$1$-cell $G \times D^1$ for
  each generator $g_{\alpha}$, attached by the left $G$-map that takes
  $\{e\} \times \partial D^1 =\{e\} \times \{0,1\}$ to $e$ and
  $g_{\alpha}$, and one $G$-$2$-cell $G \times D^2$ for each relation
  $r_{\beta}$ attached by the left $G$-map that on $\{e\} \times
  \partial D^2$ is the loop $r_{\beta}$ at $e$.  The orbit space under
  the left action of $H < G$ on ${X}_G(\{e\} \backslash G)$ is the Cayley
  complex for the orbit space $H \backslash G$: $H \backslash
  {X}_{\{e\} \backslash G} = {X}_{H \backslash G}$. In particular,
  \begin{equation*}
    G \backslash
  {X}_{\{e\} \backslash G} = {X}_{G \backslash G} = D^0 \cup
  \coprod_{\{g_{\alpha}\}} D^1 \cup \coprod_{\{r_{\beta}\}} D^2 = X_G
  \end{equation*}
  is a point $\{Ge\}$ with an arrow $Ge \xrightarrow{g_{\alpha}} Ge$
  for each generator $g_{\alpha}$ and with one $2$-cell attached along
  the loop $Ge \xrightarrow{r_{\beta}} Ge$ for each relation
  $r_{\beta}$. 
  

  It is very instructive to do a few examples. See \cite{diestel} for
  information about graph theory.

  \begin{exmp}[Cayley complexes for cyclic groups]
    For the infinite cyclic group $G=C_\infty=\langle a \rangle$,
    $X_{\{e\} \backslash G} = \Z \cup (\Z \times D^1) = \R$ and $X_{G
      \backslash G} = G \backslash \R = S^1$.  For the cyclic group
    $G=C_2=\langle a \mid a^2 \rangle$ of order $2$, $X_{\{e\}
      \backslash G} = (C_2 \times D^0) \cup (C_2 \times D^1) \cup (C_2
    \times D^2) = S^2$ and $X_{G \backslash G} = G \backslash S^2 = \R
    P^2$. For the cyclic group $G=C_m=\langle a \mid a^m \rangle$ of
    order $m$, $X_{\{e\} \backslash G}$ is a circle with $C_m \times
    D^2$ attached and and $X_{G \backslash G} = G \backslash X_{\{e\}
      \backslash G}$ is the mapping cone for $S^1 \xrightarrow{m}
    S^1$.
  \end{exmp}

  \begin{exmp}[Cayley graphs for  $F_2$-sets] 
    Let $G=\gen{a,b}=\Z \amalg \Z$ be a free group $F_2$ on two generators.
    Then $X_{G \backslash G} = S^1 \vee S^1$ and $\mathcal{O}_G =
    \mathrm{Cov}_0(S^1 \vee S^1)$.  Since there are no relations,
    Cayley complexes for right $G$-sets are Cayley graphs. In
    particular, $X_{\{e\} \backslash G} = (G \times D^0) \cup (G
    \times D^1 \coprod G \times D^1)$ is the $G$-graph
       \begin{equation*}
      \xymatrix@C=15pt@R=15pt{
                        & gb \\
         ga^{-1} \ar@{-}[r] & g \ar@{-}[u] \ar@{-}[r] & ga \\
                        & gb \ar@{-}[u] }
    \end{equation*}
    with vertex set $G$ and two edges from $g$ to $ga$ and $gb$ for
    every vertex $g \in G$, and $X_{G \backslash G} = G \backslash
    X_{\{e\} \backslash G}$ is the graph, $S^1 \vee S^1$, with one vertex $G
    \backslash G$ and two edges. In general, for any subgroup $H$ of
    $G$, the Cayley graph, $X_{H \backslash G}$, for $H \backslash G$
    is the covering space of $X_{G \backslash G} =S^1 \vee S^1$
    characterized by any of these three properties:
    \begin{itemize}
    \item $X_{H \backslash G}$ is the Cayley table for $H \backslash
      G \times G \to H \backslash G$ relative to the generators $a$ and $b$
    \item the fibre of $X_{H \backslash G} \to X_{G \backslash G}$ is
      the right $G$-set $H \backslash G$
    \item the image of the mono\m\ $\pi_1(X_{H \backslash G}) \to
      \pi_1(X_{G \backslash G}) = G$ is conjugate to $H$
    \end{itemize}  
    Here are some examples:
    \begin{itemize}
    \item If $H=\gen{a²,ab,b²}$ then the Cayley table and the Cayley
      graph of the $G$-set
      $H \backslash G=\{He,Ha\}$ are
  \begin{equation*}
    \text{
        \begin{tabular}{l|c|c}
         {} & $a$ & $b$ \\ \hline
         $He$ & $Ha$ & $Ha$ \\
         $Ha$ & $He$ & $He$
        \end{tabular}} \qquad \qquad
    \xymatrix{
      He \ar@/^1pc/[r]_a  \ar@/_1pc/[r]^b &
      Ha \ar@/_2pc/[l]^b   \ar@/^2pc/[l]_a }
  \end{equation*}
  because $He \xrightarrow{a} Ha$, $He \xrightarrow{b} Hb =
  Habb^{-2}b=Ha$, $Ha \xrightarrow{a} Haa=He$, and $Ha \xrightarrow{b}
  Hab =He$. The subgroup $H$ is normal since it has index two. Note
  that $H$ is free of rank $3$.
%
\item If $H = [G,G]$ is the commutator subgroup of $G$ then the Cayley
  graph gives a tiling of the the plane by squares with edges
  labeled $aba^{-1}b^{-1}$.
\item If $H=G^2$ is the smallest subgroup containing all squares in
  $G$, the right cosets are $H \backslash G = \{He,Ha,Hb,Hab\}$ and
  the Cayley graph is the graph of the Cayley table.
\item If $H$ is the smallest normal subgroup containing $a^3$, $b^3$,
  and $(ab)^3$, then the Cayley graph gives a tiling of the plane
  by hexagons, with edges $ababab$, and triangles with edges $aaa$
  or $bbb$. Observe that $Hxa^3y = Hxa^{-3}x^{-1}xa^3y=Hxy$.
\end{itemize}
    It is, in general, a difficult problem to enumerate the cosets of
    $H$ in $G$. 
     \end{exmp}
 
     \begin{exc} \cite[(3) p 58]{hatcher} Let $G=<a,b>$ be the free
       group on two generators and $H=\gen{a^2, b^2, aba^{-1},
         bab^{-1}}$. Draw the Cayley graph for $H \backslash G$ with
       the help of the information provided by this magma session:
\begin{verbatim}
> G<a,b>:=FreeGroup(2);
>  H:=sub<G|a^2, b^2, a*b*a^-1, b*a*b^-1>;
> Index(G,H);
3
> T,f:=RightTransversal(G,H);
> T;
{@ Id(G), a, b @} //The vertices of the Cayley graph
> E:={@ <v,(v*a)@f,(v*b)@f> : v in T @};
> E;
{@ <Id(G), a, b>, <a, Id(G), a>, <b, b, Id(G)> @} //The edges 
> 
\end{verbatim}
  \end{exc}

  \begin{exc}
    Let $G$ be a free group of finite rank and $H$ a subgroup of $G$.
    Show that $H$ is free and that $\order{G \colon
      H}(\mathrm{rk}(G)-1)= {\mathrm{rk}(H)-1}$. (This exercise is
    most easily solved by using the Euler characteristic.)
  \end{exc}

  \begin{exmp}
    When $G=C_m= \gen{g \mid g^m}$ is the cyclic group of order $m>0$,
    the Cayley complex
    \begin{equation*}
    {X}_G(\{e\} \backslash G) = 
   (G \times D^0) \cup (G \times D^1) \cup (G \times D^2)     
    \end{equation*}
    is the universal covering space of the mapping cone for the degree
    $m$ map on the circle. It is the left $G$-CW-complex consisting of a
    circle with $m$ $2$-discs attached.  (When $m=2$, this is the
    $2$-sphere which is the universal covering space of the mapping
    cone $\R P^2$ for the degree $2$-map of the circle.)  What is the
    covering space action of $G$ on ${X}_G(\{e\} \backslash G)$?
  \end{exmp}

  \begin{exmp}\label{exmp:RP2astRP2}
    \cite[Example 1.48, Exercise 1.3.14]{hatcher} Let $G=\gen{a,b \mid
      a^2,b^2}=\Z/2 \amalg \Z/2 \stackrel{\ref{exmp:amalg}}{=} \Z
    \rtimes \Z/2$ be the free product of $\Z/2$ with itself. Then
    $X_{G \backslash G} = \R P^2 \vee \R P^2$ and $\mathrm{Cov}_0(\R
    P^2 \vee \R P^2)=\mathcal{O}_{C_2 \amalg C_2}$. The total space
    $X_{\{e\} \backslash G}$ of its universal covering space
    ${X}_{\{e\} \backslash G} \to X_{G \backslash G}$ is an infinite
    string of $S^2$s. Indeed, the $0$-skeleton is $G$, the
    $1$-skeleton obtained by attaching two $1$-discs to each $0$-cell,
    is
    \begin{equation*}
      \xymatrix@1{
    {\cdots} &
     ba \ar@{-}@/_/[r]_a \ar@{-}@/^/[r]^a & 
     b \ar@{-}@/_/[r]_b \ar@{-}@/^/[r]^b &
     e \ar@{-}@/_/[r]_a \ar@{-}@/^/[r]^a &
     a \ar@{-}@/_/[r]_b \ar@{-}@/^/[r]^b &
     ab &
     {\cdots} } 
    \end{equation*}
    and the $2$-skeleton is obtained by attaching two $2$-discs at
    each $0$-cell along the maps $a^2$ and $b^2$. The left action of
    $a \in G$ which swaps $e \leftrightarrow a$, $b \leftrightarrow
    ab$, etc is the antipodal map on the sphere containing $e$ and $a$.
    
    The subgroup $H=\gen{(ab)^3}=3\Z \subset \Z$ is normal in $G$ so
    that the orbit set 
    \begin{equation*}
    H \backslash G=\{He,Ha,Hb,Hab,Hba,Haba\} = 3\Z
    \backslash \Z \rtimes 2\Z \backslash \Z     
    \end{equation*}
    is actually a group; it is the dihedral group of order $6$,
    isomorphic to $\Sigma_3$.  The quotient space $H \backslash
    {X}_{\{e\} \backslash G} = {X}_{H \backslash G}$ is a necklace of
    six $S^2$s formed from the $1$-skeleton
  \begin{equation*}
  \xymatrix{
   & He \ar@/^1pc/[dr]^a \ar@/_1pc/[dl]_b \ar@{<-->}[dr] \\
   Hb \ar@/_1pc/[ur]_b \ar@/_1pc/[d]_a \ar@{<-->}[drr]  && 
   Ha \ar@/^1pc/[ul]^a \ar@/_1pc/[d]_b \\
   Hba \ar@/_1pc/[u]_a \ar@/_1pc/[dr]_b \ar@{<-->}[dr] &&
   Hab \ar@/_1pc/[u]_b \ar@/_1pc/[dl]_a \\
   & Haba \ar@/_1pc/[ul]_b  \ar@/_1pc/[ur]_a }
   \end{equation*}
   by attaching $2$-discs at each vertex along the loops $a^2$ and
   $b^2$. The fundamental group of $H \backslash \widetilde{X}_G$ is
   $H$ and the deck transformation group is $H \backslash N_G(H) = H
   \backslash G$ since $H$ is normal.  The dashed arrows show the
   covering space left action by $Ha \in H \backslash G$; the orbit
   space for this action is the Cayley complex of the next example.
   The element $Hab \in H \backslash G$ acts by rotating the graph two
   places in clockwise direction.

  For another example,
   take $H=\gen{(ab)^3,a}=3\Z \rtimes \Z/2$; $H$ is not normal for
   $N_G(H)=H$,  $H
   \backslash G=\{He,Hb,Hba\}$ has $3$ elements, and
\begin{equation*}
  \xymatrix{
  He \ar@/^1pc/[r]^b \ar@(ul,dl)_a &
  Hb \ar@/^1pc/[l]^b \ar@/^1pc/[r]^a &
  Hba  \ar@/^1pc/[l]^a \ar@(ur,dr)^b }
\end{equation*}
is the Cayley graph for $H \backslash G$. The Cayley complex, obtained
by attaching six $2$-discs along the maps $a^2$ and $b^2$ at each
vertex, is $\R P^2$, $S^2$, $S^2$, $\R P²$ on a string as shown above.
The deck transformation group $H \backslash N_G(H)= H \backslash H $
is trivial.
\end{exmp}

\begin{figure}\label{fig1}
\centering
\includegraphics[width=5cm,height=2cm]{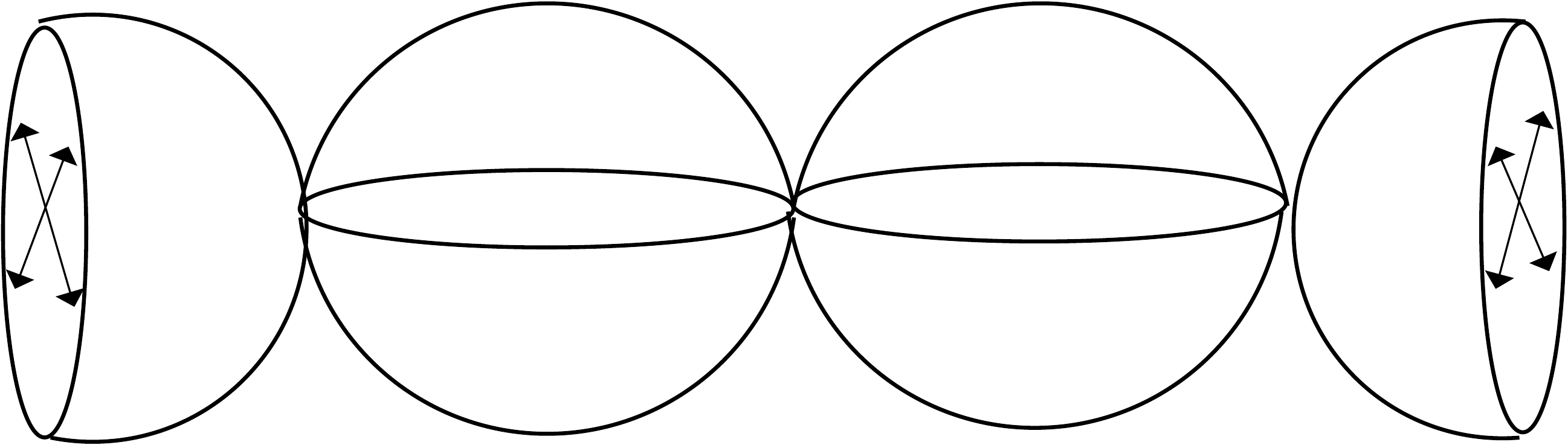}
\caption{$H \backslash {X}_{\{e\} \backslash G}$ for $H=\gen{(ab)^3,a} \leq
 \gen{a,b \mid a^2,b^2}=G$}
\end{figure} 

\begin{exmp} Let $G=\Z/2 \amalg \Z/3  \stackrel{\ref{exmp:amalg}}{=}
  \mathrm{PSL}(2,\Z)$ be the free product of a cyclic group of order
  two and a cyclic group of order three. This graph
\begin{equation*}
  \xymatrix{
   b \ar@{-}@/^/[dr]^b \ar@{-}@/_2pc/[dd]_b \\
   & e \ar@{-}@/^1pc/[r]^a \ar@{-}@/_1pc/[r]_a & a \\
   b^2 \ar@{-}@/_/[ur]_b } 
\end{equation*}
is the beginning of the Cayley complex for $G$. Describe the left
$G$-CW-complex ${X}_G(G)$!
\end{exmp}

\setcounter{subsection}{\value{thm}}
\subsection{Normal covering maps}
\label{sec:normal}\add

Let \func{p}{Y}{X} be a covering map between path connected spaces. 
\begin{defn}
  The covering map \func{p}{Y}{X} is normal if the group
  $\mathrm{Cov}(X)(Y,Y)$ of deck transformations acts transitively on
  the fibre $p^{-1}(x)$ over some point of $X$.
\end{defn}

If the action is transitive at some point, then it is transitive at
all points. Why are these covering maps called normal covering maps?  

\begin{cor}
  Let $X$ be a path connected, locally path connected and semi-locally
  simply connected space and \func pYX a covering map with $Y$ path
  connected. Then
  \begin{equation*}
    \text{The covering map \func pYX is normal\ } \iff 
    \text{The subgroup $\pi_1(Y)$ is normal in $\pi_1(X)$} 
  \end{equation*}
\end{cor}
\begin{proof}
  The action $ H \backslash N_G(H) \times H \backslash G
  \to H \backslash G$ of the group of covering maps on the fibre is
  transitive iff and only if $H$ is normal in $G$.
\end{proof}

All double covering maps are normal since all index two subgroups are
normal.

\setcounter{subsection}{\value{thm}}
\subsection{Sections in covering maps}
\label{sec:sections}\add
A {\em section\/} of a covering \func{p}{E}{X} is a (\co ) map
\func{s}{X}{E} such that $s(x)$ lies above $x$, $ps(x)=x$, for all $x
\in X$. In other words, a section is a lift of the identity map of the
base space. Each section traces out a copy of the base space in the
total space (and that is why it is called a section).

\begin{lemma}\label{eq:fixform}
  Let \func{p}{E}{X} be a covering space over a connected, locally
  path connected and semi-locally simply connected base space
  $X$. Then the evaluation map $s \to s(x)$
  \begin{equation*}
    \xymatrix@1{
      {\{\text{sections of \func{p}{E}{X}} \}} \ar[r] &  
         p^{-1}(x)^{\pi_1(X,x)} }
  \end{equation*}
  is a bijection.
\end{lemma}
\begin{proof}
  Since $X$ is connected, sections are determined by their value at a
  single point (\ref{thm:lift}), so the map is injective. It is also
  surjective because any $\pi_1(X,x)$-invariant point corresponds
  (under the classification of covering spaces over $X$) to the
  trivial covering map $X \to X$ which obviously has a section.
\end{proof}

In fact, $E$ contains the trivial covering $ p^{-1}(x)^{\pi_1(X,x)}
\times X$ as a subcovering.

If either of $Y_1 \to X$ or $Y_2 \to X$ is normal, then 
\begin{equation*}
  \mathrm{Cov}(X)(Y_1,Y_2)= 
  \begin{cases}
    \pi_1(Y_2) \backslash \pi_1(X) & \pi_1(Y_2) \subset \pi_1(Y_1) \\
    \emptyset & \text{otherwise}
  \end{cases}
\end{equation*}
for the transporter $N_{\pi_1(X)}(\pi_1(Y_1),\pi_1(Y_2))$ equals
$\pi_1(X)$ if $\pi_1(Y_1) \subset \pi_1(Y_2)$ and $\emptyset$ otherwise.



\section{Universal covering spaces of topological groups}
\label{sec:covtopgroup}

Suppose that $G$ is a connected, locally path connected, and
semi-locally simply connected topological group (for instance, a
connected Lie group) and let $G\gen{1}$ be the universal covering
space (\ref{defn:X1}) of $G$. We can use the group multiplication in
$G$ to define a multiplication in $G\gen{1}$ simply by letting the
product $[\gamma] \cdot [\eta]$ of two homotopy classes of paths
$[\gamma], [\eta] \in G\gen{1}$ equal the homotopy class $[\gamma\cdot
\eta] \in G\gen{1}$ of the product path
$(\gamma\cdot\eta)(t)=\gamma(t)\cdot\eta(t)$ whose value at any time
$t$ is the product of the values $\gamma(t)\in G$ and $\eta(t) \in G$.

\begin{lemma}
  $G\gen{1}$ is a topological group and $G\gen{1} \to G$ is a \m\ of
  topological groups whose kernel is the subgroup $\{[\omega] \mid
  \omega(0)=\omega(1)\} = \pi_1(G,e)$ of homotopy classes of loops
  based at the unit $e \in G$.
\end{lemma}

The set $\pi_1(G,e)$ is here equipped with the group structure it
inherits from $G\gen{1}$ where multiplication of paths is induced from
group multiplication in $G$.  However, we have also defined a group
structure on $\pi_1(G,e)$ using composition of loops. It turns out
that these two structures are identical.

\begin{lemma}\label{lemma:twogroupstr}
  Let $\omega_1$ and $\omega_2$ be two loops in $G$ based at the unit
  element $e$. Then the loops $\omega_1 \cdot \omega_2$ (group
  multiplication) and $\omega_1\omega_2$ (loop composition) are
  homotopic loops.
\end{lemma}
\begin{proof}
  There is a homotopy commutative diagram
  \begin{equation*}
    \xymatrix@C=45pt{
      S^1 \ar[dr]_{\Delta} \ar[r] & 
      S^1 \vee S^1 \ar@{^(->}[d] \ar[r]^{\omega_1 \vee \omega_2} &
      G \vee G \ar@{^(->}[d] \ar[dr]^{\nabla} \\
      & 
      S^1 \times S^1 \ar[r]_{\omega_1\times\omega_2} &
      G \times G \ar[r] &
      G }
  \end{equation*}
  where $\Delta$ is the diagonal and $\nabla$ the folding map. The
  loop defined by the top edge from $S^1$ to $G$ is the composite loop
  $\omega_1\omega_2$ and the loop defined by the bottom edge is the
  product loop $\omega_1\cdot\omega_2$.
\end{proof}

One can also show that in this situation $\pi_1(G,e)$ must be
abelian.

Let $\Ha=\R 1 \oplus \R i \oplus \R j \oplus \R k$ be the quaternion
algebra where the rules $i^2=j^2=k^2=-1$, $ij=k=-ji$,
$jk=i=-kj$, $ki=j=-ki$ define the multiplication. Let $\Symp(1)$
denote the topological group of quaternions of norm $1$.

$\Symp(1)$ acts in a norm preserving way on the real vector space
$\Ha=\R^4$ by the rule $\alpha \cdot v = \alpha v \alpha^{-1}$ for all
$\alpha \in \Symp(1)$ and $v \in \R^4 = \Ha$. This give a homo\m\ 
\func{\pi}{\Symp(1)}{\SO(4)}.  Since $\R 1$ is invariant under this
action, it takes $\R^{\perp} = \R i \oplus \R j \oplus \R k = \R^3$ to
itself, so there is also a group homo\m\ \func{\pi}{\Symp(1)}{\SO(3)}
\cite[I.6.18, p 88]{brockerdieck}. The kernel is $\R \cap \Symp(1) =
\{\pm 1\}$. Convince yourself that $\pi$ is surjective (see the
computation below and recall that an element of $\SO(3)$ is a rotation
around a fixed line), so that \func{\pi}{\Symp(1)}{\SO(3)= \{\pm 1\}
  \backslash \Symp(1)} is a double covering space.  

Let
\begin{equation*}
  R(\theta)=
  \begin{pmatrix}
    \cos\theta & - \sin\theta \\ \sin\theta & \phantom{-}\cos\theta
  \end{pmatrix}
\end{equation*}
be the matrix for rotation through angle $\theta$.

\begin{lemma}
  The map \func{\pi}{\Symp(1)}{\SO(3)} is the universal covering map
  of $\SO(3)$. The fundamental group $\pi_1(\SO(3),E)=\{\pm 1\}$ is
  generated by the loop 
  \begin{equation*}
  \omega(t) = \pi\alpha(t) = \begin{pmatrix}
    R(2\pi t) & 0 \\ 0 & 1
  \end{pmatrix}, \qquad 0 \leq t \leq 1,
\end{equation*}
\end{lemma}
\begin{proof}
  The topological space $\Symp(1)=S^3$ is simply connected, so
  $\Symp(1) \to \SO(3)$ is the universal covering space of $\SO(3)$.
  (We have seen this double covering before: It is the double covering
  $S^3 \to \R P^3$.)

The fundamental group $\pi_1(\SO(3),E)=C_2$ is generated by the image
loop $\omega(t)=\pi\alpha(t)$ of a path $\alpha(t)$ in $\Symp(1)$ from
$+1$ to $-1$. If we take
\begin{equation*}
  \alpha(t)= \cos(\pi t) +  \sin(\pi t)k, \qquad 0 \leq t \leq 1,
\end{equation*}
then the image in $\SO(3)$ is the loop
\begin{equation*}
  \omega(t) = \pi\alpha(t) = \begin{pmatrix}
    R(2\pi t) & 0 \\ 0 & 1
  \end{pmatrix}, \qquad 0 \leq t \leq 1,
\end{equation*}
This follows from the computation
\begin{multline*}
  \alpha(t)i\alpha(t)^{-1} = (\cos(\pi t) + k \sin(\pi t))i( (\cos(\pi
  t) - k \sin(\pi t)) \\ = \cos^2(\pi t)i + \cos(\pi t)\sin(\pi t)j + 
  \cos(\pi t)\sin(\pi t)j - \sin^2(\pi t)i 
  = \cos(2\pi t)i + \sin(2\pi t)j
\end{multline*}
and similarly for  $\alpha(t)j\alpha(t)^{-1}=-\sin(2\pi t)i + \cos(2\pi
t)j$ and
$\alpha(t)k\alpha(t)^{-1}=k$. 
\end{proof}

 
It is also known that the inclusion $\SO(3) \to \SO(n)$ induces an
iso\m\ on $\pi_1$ for $n\geq 3$. We conclude that the fundamental
group $\pi_1(\SO(n),E)$ has order two for all $n \geq 3$ and that it
is generated by the loop $\omega(t)$ in $\SO(n)$. Thus the topological
groups $\SO(n)$, $n\geq 3$, have universal {\em double\/} covering
spaces that are topological groups.

\begin{defn}
  For $n \geq 3$, $\Spin(n) = \SO(n)\gen{1}$ is the universal covering
  space of $\SO(n)$ and \func{\pi}{\Spin(n)}{\SO(n)} is the universal
  covering map.
\end{defn}

The elements of $\Spin(n)$ are homotopy classes of paths in $\SO(n)$
starting at $E$ and, in particular, $\Spin(3)=\Symp(1)$. The kernel of
the homo\m\ $\pi \colon \Spin(n) \to \SO(n)$ is $\{e,z\}$ where $e$ is
the unit element and $z=[\omega]$ is the homotopy class of the loop
$\omega$.

\begin{prop}
  The center of $\Spin(n)$ is 
  \begin{equation*}
    Z(\Spin(n)) = 
    \begin{cases}
      C_2 =\{e,z\} & \text{$n$ odd} \\
      C_2 \times C_2=\{e,x\} \times \{e,z\} & n \equiv 0 \bmod 4 \\
      C_4 = \{e,x,x^2,x^3\} & n \equiv 2 \bmod 4 
    \end{cases}
  \end{equation*}
  for $n \geq 3$. 
\end{prop}
\begin{proof}
  From Lie group theory we know that the center of $\Spin(n)$ is the
  inverse image of the center of $\SO(n)$. Thus the center of
  $\Spin(n)$ has order $2$ when $n$ is odd and order $4$ when $n$ is
  even.
  
  Suppose that $n=2m$ is even. Then $Z(\Spin(2m))=\{e,z,x,zx\}$ where
  $x=[\eta]$ is the homotopy class of the path
\begin{equation*}
  \eta(t) = \diag(R(\pi t), \ldots , R(\pi t))
\end{equation*}
from $E$ to $-E$. Note  that $x^2$ is (\ref{lemma:twogroupstr})
represented by the loop
\begin{equation*}
  \eta(t)^2 = \diag(R(2\pi t), \ldots , R(2\pi t))
\end{equation*}
Conjugation with a permutation matrix from $\SO(2n)$ takes 
\begin{equation*}
  \omega(t)=\diag(R(2\pi t),E,E,\ldots,E) \text{\ to\ }
   \diag(E,R(2\pi t),E,\ldots,E)
\end{equation*}
and since inner auto\m s are based homotopic to identity maps, 
both the above loops represent the generating loop $\omega$. It
follows that 
\begin{equation*}
  x^2 =[\eta(t)^2] = [\omega(t)^m] = z^m= 
  \begin{cases}
    e & \text{$m$ even} \\
    z  & \text{$m$ odd}
  \end{cases}
\end{equation*}
Thus $Z(\Spin(2m))=\{z \} \times \{ x\} = C_2 \times C_2$
if $m$ is even and $Z(\Spin(2m))=\{x\} = C_4$ if $m$ is odd. 
\end{proof}

What is the fundamental group $\pi_1(\PSO(2n))$ of the topological
group $\PSO(2n)=\SO(2n)/\gen{-E}$? 

When will two diagonal matrices in $\SO(n)$ commute in $\Spin(n)$? Let
$D=\{\diag(\pm 1, \ldots, \pm 1\}$ be the abelian subgroup of diagonal
matrices in $\SO(n)$. By computing commutators and squares in
$\Spin(n)$ we obtain functions 
\begin{equation*}
  \func{[\ ,\ ]}{D \times D}{\{e,z\}}, \qquad
  \func{q}{D}{\{e,z\}}
\end{equation*}
given by $q(d)=(\bar d)^2$ and $[d_1,d_2] = [\bar d_1, \bar d_2]$ where
$\pi(\bar d)=d$, $\pi(\bar d_1)=d_1$, $\pi(\bar d_2)=d_2$. They are
related by formula
\begin{equation*}
  q(d_1+d_2)=q(d_1)+q(d_2)+[d_1,d_2]
\end{equation*}
which says that $[\ , \ ]$ records the deviation from $q$ being a group
homo\m\ (using additive notation here). It suffices to compute $q$ in
order to answer the question about commutativity relations.

\begin{prop}
  $q(d)=e$ iff the number of negative entries in the diagonal matrix
  $d \in D$ is divisible by $4$. $[d_1,d_2]=e$ iff the number of
  entries that are negative in both $d_1$ and $d_2$ is even.
\end{prop}
\begin{proof}
  Note that two elements of $D$ are conjugate iff they have the same
  number of negative entries. Use permutation matrices and, if
  necessary, the matrix $\diag(-1,1,\ldots,1)$. Consider for instance 
  \begin{equation*}
    d_1 = \diag(-1,-1,1,\ldots,1), \qquad
    d_2=\diag(-1,-1,-1,-1,1,\ldots,1) 
  \end{equation*}
  with two, respectively four, negative entries. The paths
  \begin{equation*}
    \bar{d}_1(t)=\diag(R(\pi t),1,\ldots,1), \qquad
    \bar{d}_2(t)=\diag(R(\pi t),R(\pi t),1,\ldots,1)
  \end{equation*}
  represent lifts of $d_1$ and $d_2$ to $\Spin(n)$. Then
  $({\bar{d}_1})^2=z=q(d_1)$ and $({\bar{d}_2})^2=e=q(d_2)$.
  Computations like these prove the formula for $q$ and the formula
  for $[\ , \ ]$ follows. The number of negative entries in $d_1+d_2$
  is the number of negative entries in $d_1$ plus the number of
  negative entries in $d_2$ minus twice the number of entries that are
  negative in both $d_1$ and $d_2$.
\end{proof}



\begin{exc}
  Let $\overline{D}_{n} \subset \Spin(n)$ be the inverse image of $D
  \subset \SO(n)$.  How many elements of order $4$ are there in
  $\overline{D}_{n}$? Can you identify the group $\overline{D}_n$?
  Show that there is a homo\m\ $\SU(m) \to \Spin(2m)$. When $m$ is
  even, what is the image of $-E \in \SU(m)$? What is the image of the
  center of $\SU(m)$? Describe the covering spaces of $\U(n)$.
\end{exc}






The inclusions $\SO(n) \subset \SO(n+1)$, $n>2$, and $\SO(m) \times
\SO(n) \subset \SO(m+n)$, $m,n>2$, of special orthogonal groups lift
to inclusions
\begin{equation*}
  \xymatrix{
   {\Spin(n)} \ar@{^(->}[r] \ar[d] & {\Spin(n+1)} \ar[d] &
   {\Spin(m) \times_{\gen{(z_1,z_2)}}\Spin(n)}
      \ar@{^(->}[r] \ar[d] & {\Spin(m+n)} \ar[d] \\
   {\SO(n)} \ar@{^(->}[r] & {\SO(n+1)} &
   {\SO(m) \times \SO(n)}  \ar@{^(->}[r] & {\SO(m+n)} }
\end{equation*}
of double coverings. (Here, $\Spin(m) \times_{\gen{(z_1,z_2)}}
  \Spin(n)$ stands for $\gen{(z_1,z_2)} \backslash (\Spin(m) \times
    \Spin(n))$).

The inclusion $\U(n) \subset \SO(2n)$, that comes from the
identification $\C^n=\R^{2n}$, lifts to an inclusion of double
covering spaces as shown in the following diagrams.
\begin{equation*}
  \xymatrix{
    { \SU(n) \times_{C_k} \U(1)} \ar[d]_{(A,z) \to (A,z)}
    \ar[r] & {\Spin(2n)} \ar[d] & 
    {\SU(n) \times_{C_n'} \U(1)} \ar[d]_{(A,z) \to (A,z^2)} 
    \ar[r] & {\Spin(2n)} \ar[d] \\
    {\SU(n) \times_{C_n} \U(1) = \U(n)} \ar@{^(->}[r] &
    {\SO(2n)} &
    {\SU(n) \times_{C_n} \U(1) = \U(n)} \ar@{^(->}[r] &
    {\SO(2n)} }
\end{equation*}
To the left, $n=2k$ is even, and to the right, $n=2k+1$ is odd;
$C_n=\{(\zeta E,\zeta^{-1}) \mid \zeta^n=1\}$ and $C_n'=\{(\zeta
E,\zeta^k) \mid \zeta^n=1\}$ are cyclic groups of order $n$ and
$C_k=\{(\zeta E,\zeta^{-1}) \mid \zeta^k=1\} \subset C_{2k}=C_n$ is
cyclic of order $k$. The iso\m\ $\SU(n) \times_{C_n} \U(1) \to \U(n)$
takes $(A,z)$ to $zA$.  When $n$ is divisible by $4$, $z=(-E,-1)$ and
$x=(E,-1)$ have order two; when $n$ is even and not divisible by $4$,
$x=(E,i)$ has order four and $x^2=(E,-1)=z$. This explains the
computation of the center of $\Spin(2n)$. (Is the group in the upper
left corner of the right diagram isomorphic to $\U(n)$? See
\cite{baum} for more information.)

There is a double covering map $\textmd{pin}(n) \to \Or(n)$ obtained
as the pullback of $\Spin(2n) \to \SO(2n)$ along the inclusion homo\m\ 
$\Or(n) \subset \SO(2n)$.

\begin{exmp}
 The inclusion $\U(2) \subset \SO(4)$ lifts to an inclusion $\SU(2)
 \times \U(1) \subset \Spin(4)$. 
 Let $G_{16} \subset \SU(2) \times \U(1) \subset \Spin(4)$ be the group
 \begin{equation*}
   G_{16} = \left\langle \left( 
       \begin{pmatrix}
         -i&0\\0&i
       \end{pmatrix},i \right), \left( 
       \begin{pmatrix}
         0&i\\i&0
       \end{pmatrix},-i \right) \right\rangle
 \end{equation*}
 $G_{16}$ has order $16$, center $Z(G_{16})=\{x,z\}=C_2 \times
 C_2=Z(\Spin(4))$, and derived group $[G_{16},G_{16}]=\{xz\}=C_2$. Its
 image under the covering maps
  \begin{equation*}
    \xymatrix{
      & {\Spin(4)/\gen{z}=\SO(4)} \ar[dr] \\
     {\Spin(4)} \ar[ur] \ar[r] \ar[dr] & {\Spin(4)/\gen{xz}} \ar[r] &
     {\Spin(4)/\gen{x,z}=\PSO(4)} \\
      & {\Spin(4)/\gen{x} =\operatorname{SSpin}(4) } \ar[ur] }
  \end{equation*}
  is dihedral $D_8$ in $\SO(4)$, abelian $C_4 \times C_2$ in
  $\Spin(4)/\gen{xz}$, quaternion $Q_8$ in the semi-spin group
  $\operatorname{SSpin}(4)=\Spin(4)/\gen{x}$, and elementary abelian
  $C_2 \times C_2$ in $\PSO(4)$. (All proper subgroups of $G_{16}$ are
  abelian but itself and some of its quotient groups are nonabelian.)
\end{exmp}

\begin{exmp}
  There exists a covering space homo\m s of topological groups
  \begin{equation*}
    \U(1) \times \SU(n) \to \U(n) \colon (z,A) \to zA
  \end{equation*}
  with kernel $C_n=\{(z,z^{-1}E) \mid z^n=1
  \}=\gen{(\zeta,\zeta^{-1}E)}$ where $\zeta_n=e^{2\pi i/n}$. The
  universal covering space homo\m\ is $\R \times \SU(n) \to \U(n)
  \colon (t,A) \to \zeta_n^tA$ with kernel
  $C_{\infty}=\gen{(1,\zeta_n^{-1}E)}$. Any covering space of $\U(n)$
  is of the form $\gen{(k,\zeta_n^{-k})} \backslash (\R \times
  \SU(n))$ for some integer $k \geq 0$. 
  
  Similarly, let $S(\U(m) \times \U(n))$ denote the closed topological
  subgroup $(\U(m) \times \U(n)) \cap \SU(m+n)$ of $\U(m+n)$.
  There exists a covering space homo\m s of topological groups
  \begin{equation*}
    \U(1) \times \U(m) \times \U(n) \to S(\U(m) \times \U(n)) 
    \colon (z,A,B) \to \diag(z^{n_1}A,z^{-m_1}B)
  \end{equation*}
  with kernel $C_{\lcm(m,n)} = \{(z,z^{-n_1}E,z^{m_1}E) \mid
  z^{\lcm(m,n)}=1 \} = \gen{(\zeta_{\lcm(m,n)},\zeta_m^{-1},\zeta_n)}$
  where $m_1=m/\gcd(m,n)=\lcm(m,n)/n$ and
  $n_1=n/\gcd(m,n)=\lcm(m,n)/m$. The universal covering space homo\m\ 
  of $ S(\U(m) \times \U(n))$ is $\R \times \SU(m) \times \SU(n) \to
  S(\U(m) \times \U(n)) \colon (t,A,B) \to (\zeta_m^tA,\zeta_n^tB)$
  with kernel $C_{\infty}=\gen{(1,\zeta_m^{-1},\zeta_n)}$.  Any
  covering space of $S(\U(m) \times \U(n)$ is of the form
  $\gen{(k,\zeta_m^{-k},\zeta_n^{k})} \backslash (\R \times \SU(m)
  \times \SU(n))$ for some integer $k \geq 0$.
\end{exmp} 

All finite covering spaces of $\U(n)$ are covered by $\U(1) \times
\SU(n)$. To see this, let $n$ and $k$ be integers and put
$k_1=k/\gcd(n,k)$. Then there is a commutative diagram
\begin{equation*}
  \xymatrix@C=50pt{
  {\U(1) \times \SU(n)} \ar@{..>}[r] \ar[d]_{(z,A)\to (z^{k_1},A)} &
  {\R \times_{\gen{(k,\zeta^{-k})}}\SU(n)} \ar[d]\\
  {\U(1) \times \SU(n)} \ar[r]_-{(\zeta_n^t,A) \to (t,A)} & 
  {\R \times_{\gen{(1,\zeta^{-1})}}\SU(n)=\U(n)} }
\end{equation*}
of covering space homo\m s.


\def\cprime{$'$}
\providecommand{\bysame}{\leavevmode\hbox to3em{\hrulefill}\thinspace}
\providecommand{\MR}{\relax\ifhmode\unskip\space\fi MR }
\providecommand{\MRhref}[2]{%
  \href{http://www.ams.org/mathscinet-getitem?mr=#1}{#2}
}
\providecommand{\href}[2]{#2}

\end{document}